\documentclass[aap,MSNbibl,dvips]{arximspdf}
\usepackage{graphicx}

\doi{10.1214/10-AAP720}
\volume{21}
\issue{3}
\pubyear{2011}
\firstpage{1102}
\lastpage{1135}

\makeatletter
\newcommand{\eqref}[1]{(\ref{#1})}

\newcommand{\ep}{\varepsilon}
\newcommand{\R}{\mathbb R}
\newcommand{\E}{\mathbb E}
\newcommand{\Prob}{\mathbb P}
\newcommand{\barq}{\overline{q}}

\newtheorem{theorem}[proposition]{Theorem}
\newtheorem{lemma}[proposition]{Lemma}
\newtheorem{corollary}[proposition]{Corollary}
\newtheorem{proposition}{Proposition}[section]

\newtheorem{hypothesis}[proposition]{Hypothesis}
\newproclaim{eg}{Example}[section]
\newproclaim{remark}{Remark}

\makeatother

\begin{document}
\begin{frontmatter}

\title{Recovering a
time-homogeneous stock price process
from perpetual option prices}
\runtitle{Recovering a time-homogeneous stock price process}

\begin{aug}
\author[A]{\fnms{Erik} \snm{Ekstr\"om}\corref{}\ead[label=e1]{ekstrom@math.uu.se}}
\and
\author[B]{\fnms{David} \snm{Hobson}\ead[label=e2]{d.hobson@warwick.ac.uk}}
\runauthor{E. Ekstr\"om and D. Hobson}
\affiliation{Uppsala University and University of Warwick}
\address[A]{Department of Mathematics\\ Uppsala University\\
Box 480,
SE-751 06 Uppsala\\ Sweden\\
\printead{e1}} 
\address[B]{Department of Statistics\\ University of Warwick\\
Coventry\\
CV4 7AL\\ United Kingdom\\
\printead{e2}}
\end{aug}

\received{\smonth{3} \syear{2009}}
\revised{\smonth{4} \syear{2010}}

%
\begin{abstract}
It is well known how to determine the price of perpetual American
options if the underlying stock price is a time-homogeneous diffusion.
In the present paper we consider the inverse problem, that is, given
prices of perpetual American options for different strikes, we show how
to construct a time-homogeneous stock price model which reproduces the
given option prices.
\end{abstract}

%

\end{frontmatter}
%

\section{Introduction}

In the classical Black--Scholes model, there is a one-to-one
correspondence between the price of an option and the volatility of the
underlying stock. If the volatility $\sigma$ is assumed to be given
(e.g., by estimation from historical data), then the arbitrage free
option price can be calculated using the Black--Scholes formula.
Conversely, if an option price is given, then the implied volatility
can be obtained as the unique $\sigma$ that would produce this option
price if inserted in the Black--Scholes formula. It has been well
documented that if the implied volatility is inferred from real market
data for option prices with the same maturity date but with different
strike prices, then, typically, a~nonconstant implied volatility is
obtained. Since the implied volatility often resembles a smile if
plotted against the strike price, this phenomenon is referred to as the
\textit{smile effect}. The smile effect is one indication that the
Black--Scholes assumption of normally distributed log-returns is too
simplistic.

A wealth of different stock price models have been proposed in order to
overcome the shortcomings of the standard Black--Scholes model, of
which the most popular are jump models and stochastic volatility
models. Given a model, option prices can be determined as risk-neutral
expectations. However, models are typically governed by a small number
of parameters, and only in exceptional circumstances can they be
calibrated to perfectly fit the full range of options data.\vadjust{\goodbreak}

Instead, there is a growing literature which tries to reverse the
procedure, using option prices to make inferences about the underlying
price process. At one extreme, models exist which take a price surface
as the initial value of a Markov process on a space of functions. In
this way the Heath--Jarrow--Morton~\cite{HJM} interest rate models can
be made to perfectly fit an initial term structure. Such ideas inspired
Dupire~\cite{Dup} to introduce the local volatility model which
calibrates perfectly to an initial volatility surface. For a local
volatility model, Dupire derived the PDE
\[
C_T(T,K)+rKC_K(T,K) =\tfrac{1}{2}\sigma^2(T,K) K^2 C_{KK}(T,K),
\]
where $C(T,K)$ is the European call option price, $T$ is time to
maturity and $K$ is the strike price. Solving for the (unknown) local
volatility $\sigma(T,K)$ gives a formula for the time-inhomogeneous
local volatility in terms of derivatives of the observed European call
option prices.

The local volatility model gives the unique martingale diffusion which
is consistent with observed call prices (alternative, nondiffusion
models also exist; see, e.g., Madan and Yor~\cite{MY}). The recent
literature (e.g., Schweizer and Wissel~\cite{SW}) has included
attempts to extend the theory to allow for a stochastic local
volatility surface. However, it relies on the knowledge of a double
continuum of option prices, which are smooth. In contrast,
Hobson~\cite{Hob} constructs models which are consistent with a
continuum of strikes, but at a single maturity, in which case there is
no uniqueness.

In the current article we present a method to recover a
time-homogeneous local volatility function from perpetual American
option prices. More precisely, we assume that perpetual put option
prices are observed for all different values of the strike price, and
we derive a time-homogeneous stock price process for which theoretical
option prices coincide with the observed ones.

No-arbitrage enforces some fundamental convexity and monotonicity
conditions on the put prices, and if these fail, then no model can
support the observed prices. If the observed put prices are smooth,
then we can use the theory of differential equations to determine a
diffusion process for which the theoretical perpetual put prices agree
with the observed prices, and our key contribution in this case is to
give an expression for the diffusion coefficient of the underlying
model in terms of the put prices. It turns out that this expression
uniquely determines the volatility coefficient at price levels below
the current stock price, but there is some freedom in the choice of the
volatility function above the current stock price level. The key idea
is to construct a dual function to the perpetual put price, and then
the diffusion coefficient can be easily found
by
taking derivatives of this dual.

The second contribution of this paper is to give time-homogeneous
models which are consistent with a given set of perpetual put prices,
even when those put prices are not twice differentiable or not strictly
convex in the continuation region where it is not optimal to exercise
immediately or not strictly convex in the continuation region. Again, the key is the dual function, coupled with a
change of scale and a time change. We give a construction of a
time-homogeneous process consistent with put prices, which we assume to
satisfy the no-arbitrage conditions, but which otherwise has no
regularity properties.

One should perhaps note that in reality, put prices are only given in
the market for a discrete set of strike prices. Therefore, as a first
step one needs to interpolate between the strikes. If a stock price is
modeled as the solution to a stochastic differential equation with a
continuous volatility function, then the perpetual put price exhibits
certain regularity properties with respect to the strike price.
Therefore, if one aims to recover a continuous volatility, then one has
to use an interpolation method that produces option prices exhibiting
this regularity. On the other hand, if a linear spline method is used,
then a continuous volatility cannot be recovered. This is one of the
motivations for searching for price processes which are consistent with
a general perpetual put price function (which is convex, but may be
neither strictly convex nor smooth).

While preparing this manuscript we came across a preprint by Alfonsi
and Jourdain, now published as~\cite{AJ08b}. The aim of~\cite{AJ08b},
as in this article, is to construct a time-homogeneous process which is
consistent with observed put prices. However, the method is different
and considerably less direct. Alfonsi and Jourdain~\cite{AJ08b}
construct a parallel model such that the put price function in the
original model (expressed as a function of strike) becomes a call price
function expressed as a function of the initial value of the stock.
They then solve the perpetual pricing problem for this parallel model
and, subject to solving a differential equation for the optimal
exercise boundaries in this model, give an analytic formula for the
volatility coefficient. In contrast, the approach in this paper is much
simpler and, unlike the method of Alfonsi and Jourdain, extends to the
irregular case.

\section{The forward problem}\label{forward}

Assume that the stock price process $X$ is modeled under the pricing
measure as the solution to the stochastic differential equation
\[
dX_t=r X_t\,dt + \sigma(X_t)X_t\,dW_t,\qquad X_0 = x_0.
\]
Here, the interest rate $r$ is a positive constant, the level-dependent
volatility $\sigma\dvtx(0,\infty)\to(0,\infty)$ is a given continuous
function and $W$ is a standard Brownian motion. We assume that the
stock pays no dividends, and we let zero be an absorbing barrier for
$X$. If the current stock price is $x_0$, then the price of a perpetual
put option with strike price $K>0$ is
%
\begin{equation}\label{hatP}
\hat P(K)=\sup_{\tau} \E^{x_0}[
e^{-r\tau}(K-X_\tau)^+],
\end{equation}
where the supremum is taken over random times $\tau$ that are stopping
times with respect to the filtration generated by $W$.
From the boundedness, monotonicity and convexity of the payoff, we have
the following.
\begin{proposition}\label{simpleproperties}
The function $\hat P\dvtx(0,\infty)\to[0,\infty)$
satisfies:
\begin{enumerate}[(ii)]
\item[(i)] $(K-x_0)^+\leq\hat P(K)\leq K$ for all $K$;
\item[(ii)]
$\hat P$ is nondecreasing and convex.
\end{enumerate}
\end{proposition}

\begin{eg*}
If $\sigma$ is constant, that
is, if $X$ is a geometric Brownian
motion, then
%
\begin{equation}\label{egBS1}
\hat P(K)=\cases{
\displaystyle\frac{K}{\beta+1}\bigl(\beta K/x_0(\beta+1)\bigr)^\beta, &\quad if $K<
\hat K$, \cr
K-x_0, &\quad if  $K\geq\hat K$,
}
\end{equation}
where $\beta= 2r/\sigma^2$ and $\hat K=x_0(\beta+1)/\beta$.
\end{eg*}

Intimately connected with the solution of the optimal stopping problem~\eqref{hatP} is the ordinary differential equation
%
\begin{equation}\label{ODE}
 \tfrac{1}{2}\sigma(x)^2 x^2u_{xx}+rxu_x-ru=0
\end{equation}
for $x>0$. This equation has two linearly independent positive
solutions which are uniquely determined (up to multiplication with
positive constants) if one requires one of them to be increasing and
the other decreasing; see, for example, Borodin and
Salminen~\cite{BS}, page 18. We denote the increasing solution by
$\hat\psi$ and the decreasing one by $\hat\varphi$. In the current
setting, $\hat\psi$ and $\hat\varphi$ are given by
\[
\hat\psi(x)=Cx
\]
and
%
\begin{equation}\label{varphi}
\hat\varphi(x)=Dx\int_x^\infty\frac{1}{y^2}
\exp\biggl\{-\int_{x_0}^y\frac{2r}{z\sigma(z)^2}\,dz\biggr\}\,dy
\end{equation}
for some arbitrary positive constants $C$ and $D$. For simplicity, and
without loss of generality, we choose
\[
D=\biggl( x_0\int_{x_0}^\infty\frac{1}{y^2}
\exp\biggl\{-\int_{x_0}^y\frac{2r}{z\sigma(z)^2}\,dz\biggr\}\,dy\biggr)^{-1}
\]
so that $\hat\varphi(x_0)=1$.

\begin{lemma}\label{varphiconvex}
The function $\hat\varphi$ is strictly decreasing
and strictly convex.
\end{lemma}

\begin{pf}
Straightforward differentiation yields
%
\begin{equation}\label{phi-derivative}
\hat\varphi'(x)=D\!\int_{x}^\infty\!\frac{1}{y^2}
\biggl(\exp\biggl\{ -\!\int_{x_0}^y\!\frac{2r}{z\sigma(z)^2}\,dz\biggr\}-
\exp\biggl\{ -\!\int_{x_0}^x\!\frac{2r}{z\sigma(z)^2}\,dz\biggr\}
\biggr)\,dy,\hspace*{-25pt}\vadjust{\goodbreak}
\end{equation}
so $\hat\varphi'(x)<0$. Similarly,
\[
\hat\varphi''(x)= \frac{2Dr}{x^2\sigma(x)^2}
\exp\biggl\{ -\int_{x_0}^x\frac{2r}{z\sigma(z)^2}\,dz\biggr\}>0,
\]
so
$\hat\varphi$ is strictly convex.
\end{pf}

\begin{figure}

\includegraphics{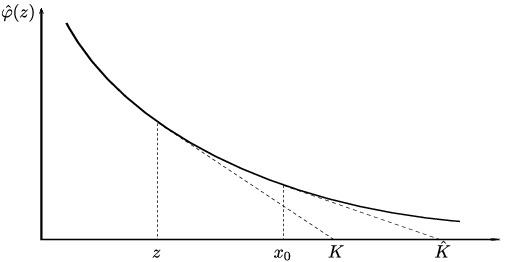}

\caption{For a given $K\leq\hat K$ the price $\hat P(K)$ is
minus the reciprocal of the slope of the tangent line to
$\hat{\varphi}$ which passes through the point
$(K,0)$.}\label{fig1}\vspace*{-3pt}
\end{figure}

It is well known that with $H_z=\inf\{t\geq0\dvtx X_t= z\}$, we
have
\begin{equation} \label{laplace}
\E^x[e^{-rH_z}]=\cases{
\hat\varphi(x)/\hat\varphi(z), &\quad if $z<x$,\vspace*{2pt}\cr
\hat\psi(x)/\hat\psi(z), &\quad if $z>x$,
}
\end{equation}
where the superindex $x$ denotes that the expected value is calculated
using $X_0=x$. [This result is easy to check by considering $e^{-r t}
\hat{\varphi}(X_t)$ and $e^{-r t} \hat{\psi}(X_t),$ which, since they
involve solutions to~\eqref{ODE}, are local martingales.] Given the
assumed time-homogeneity of the process $X$, it is natural to consider
stopping times in~\eqref{hatP} that are hitting times. Define
%
\begin{eqnarray}\label{Ptilde}
\tilde{P}(K) &:=& \sup_{z:z\leq x_0 \wedge K}
\E^{x_0}[ e^{-rH_z}(K-X_{H_z})^+]\nonumber\\
&=& \sup_{z:z\leq x_0 \wedge K}(K-z)\E^{x_0}[e^{-rH_z}] \\
&=& \sup_{z:z\leq x_0 \wedge K}\frac{K-z}{\hat
\varphi(z)},\nonumber
\end{eqnarray}
where the last equality follows from~\eqref{laplace}. Clearly
$\hat{P}(K) \geq\tilde{P}(K)$, and, of course, as we show below, there
is equality. Since the function $\hat\varphi$ is strictly convex, for
each fixed $K$ there exists a unique $z=z(K)\leq x_0$ for which the
supremum in~\eqref{Ptilde} is attained, that is,
\begin{equation}\label{PzK}
\tilde P(K)=\frac{K-z(K)}{\hat\varphi(z(K))}.
\end{equation}
Geometrically, $z=z(K)$ is the unique value (less than or equal to
$x_0$) which makes the negative slope of the line through $(K,0)$ and
$(z,\hat\varphi(z))$ as large as possible; see Figure~\ref{fig1}.

%

Define
\[
\hat K:=x_0-1/\hat\varphi^\prime(x_0).
\]
From the strict convexity of $\hat\varphi$ it follows that if $K\geq
\hat K$, then
\[
\tilde P(K)= \sup_{z:z\leq x_0}\frac{K-z}{\hat\varphi(z)}=
\frac{K-x_0}{\hat\varphi(x_0)} = K-x_0,
\]
and if $K\leq\hat K$, then
%
\begin{equation}\label{Ptilde2}
\tilde P(K)=\sup_{z:z\leq x_0}\frac{K-z}{\hat\varphi(z)}= \sup_{z}
\frac{K-z}{\hat\varphi(z)}.
\end{equation}
Moreover, for $K<\hat K$ we have $\tilde P(K)>(K-x_0)^+$.

\begin{lemma}\label{Phat}
The functions $\hat P$ and $\tilde P$ coincide, that
is,
\begin{equation} \label{Px0}
\hat{P}(K) = \sup_{z:z\leq x_0}
\frac{K-z}{\hat\varphi(z)}.\vspace*{-1pt}
\end{equation}
\end{lemma}

\begin{pf}
As noted above, we clearly have $\hat P\geq\tilde P$ since the
supremum over all stopping times is at least as large as the supremum
over first hitting times.

For the reverse implication, suppose first that $K \leq\hat{K}$. In
that case $\hat{\varphi}(z) \geq(K-z)^+/ \tilde{P}(K)$, by
\eqref{Ptilde2}. Further, $e^{-rt} \hat{\varphi}(X_t)$ is a nonnegative
local martingale and hence a supermartingale. Thus, for any stopping
time $\tau$ we have
\[
1 \geq\E^{x_0}[ e^{-r \tau} \hat{\varphi}(X_\tau) ] \geq\E^{x_0}[
e^{-r \tau} (K-X_\tau)^+/\tilde{P}(K) ] .
\]
Hence, $\tilde{P}(K) \geq
\sup_{\tau} \E^{x_0}[ e^{-r \tau} (K- X_\tau)^+] = \hat{P}(K)$.\vspace*{1pt}

Finally, let $K>\hat{K}$. It follows from the first part that $\hat
P(\hat K)=\hat K-x_0$, so Proposition~\ref{simpleproperties} implies
that $\hat P(K)=K-x_0=\tilde P(K)$, which completes the proof.\vspace*{-1pt}
\end{pf}

\begin{eg*}
If $\sigma$ is constant, that
is, if $X$ is a geometric Brownian
motion, then
\[
\hat\varphi(x)= \biggl( \frac{x_0}{x} \biggr)^{\beta},
\]
where $\beta=2r/\sigma^2$. Consequently, the put option price is given
by
\[
\hat P(K)=x_0^{-\beta}\sup_{z:z\leq x_0} (K-z)z^\beta.
\]
Straightforward differentiation shows that the supremum is attained for
\[
z=z^*:=\frac{\beta K}{\beta+1}
\]
if $z^*<x_0$, and for $z=x_0$ if $ z^*\geq x_0$. Consequently, $\hat
P(K)$ is given by~\eqref{egBS1}.\vadjust{\goodbreak}
\end{eg*}

Under our current assumptions it is not possible to rule out the case
where the diffusion $X$ hits zero in finite time, although we then
insist that zero is absorbing. Note that $X$ hits zero in finite time
if and only if $\hat\varphi(0)<\infty$, in which case we set
$\underline{K}= - \hat\varphi(0)/\hat\varphi^\prime(0)$. When
$\hat\varphi^\prime(0)$ is finite we have $\underline{K}>0$ and for
$K<\underline{K}$, $z(K)=0$ and $\hat{P}(K)=K/\hat\varphi(0)$. By the
strict concavity of~$\hat\varphi$, $\lim_{K \downarrow\underline{K}}
z(K)=0$.

\begin{proposition}
\label{Kstar} \label{smoothproperties} In addition to the properties
described in Proposi\-tion~\ref{simpleproperties}, the following
statements about the function $\hat P\dvtx[0,\infty)\to[0,\infty)$~hold:

\begin{enumerate}[(iii)]
\item[(i)] $\hat P$ satisfies $\hat P(K)>(K-x_0)^+$ for all $K\in
(0,\hat K)$ and $\hat P(K)=K-x_0$ for all $K\geq\hat K$;
\item[(ii)]
$\hat P$ is continuously differentiable on $(0,\infty)$ and twice
continuously differentiable on $(0,\infty)\setminus\{\underline
{K},\hat
K\}$;
\item[(iii)] $\hat P$ is strictly increasing on $(0,\infty)$ with
a strictly positive second derivative on $(\underline{K},\hat K)$.
\end{enumerate}
\end{proposition}

\begin{pf}
Statement (i) follows from Lemma~\ref{Phat} and the fact that (i) is
true for~$\tilde P$.

Next, consider $\underline{K}<K<\hat K$. By~\eqref{Ptilde2} we have
\[
\hat P(K)=\sup_z\frac{K-z}{\hat\varphi(z)}=
\frac{K-z(K)}{\hat\varphi(z(K))}
\]
for some $z(K)\in(0,x_0)$. Since
$z=z(K)$ maximizes the quotient $(K-z)/\hat\varphi(z)$, we have
%
\begin{equation}\label{impl}
\bigl(K-z(K)\bigr)\hat\varphi^\prime(z(K)) +\hat\varphi(z(K))=0.
\end{equation}
It follows from~\eqref{impl} and the implicit function theorem that
$z(K)$ is continuously differentiable for $\underline{K}<K<\hat K$.
Therefore, differentiating~\eqref{PzK} gives
%
\begin{eqnarray}\label{Pprime}
\hat P^\prime(K) &=&
\frac{(1-z^\prime(K))\hat\varphi(z(K))-
(K-z(K))z^\prime(K)\hat\varphi^\prime(z(K))}{(\hat\varphi(z(K)))^2}\nonumber\\ [-9pt]\\ [-9pt]
&=& \frac{1}{\hat\varphi(z(K))},\nonumber
\end{eqnarray}
where the second equality follows from~\eqref{impl}. Equation
\eqref{Pprime} shows that $\hat P^\prime(\hat K-)=
1/\hat\varphi(x_0)=1$, so $\hat P$ is $C^1$ at $\hat K,$ and, again,
when $\underline{K}>0$ we have $\hat P^\prime(\underline{K}+)=
1/\hat\varphi(0+),$ so $\hat P$ is $C^1$ also at $\underline{K}$.
Moreover, since $\hat\varphi(z)$ is $C^1$\vspace*{1pt} and $z(K)$ is $C^1$ away from
$\hat K$, it follows that $\hat P(K)$ is $C^2$ on
$(0,\infty)\setminus\{ \underline{K}, \hat K\}$. In fact, for
$\underline{K}<K<\hat K$ we have
\begin{eqnarray*}
\hat P^{\prime\prime}(K) &=&
\frac{-z^\prime(K)\hat\varphi^\prime(z(K))}
{(\hat\varphi(z(K)))^2}\\[-2pt]
&=& \frac{(\hat\varphi^\prime(z(K)))^2}
{(K-z(K))(\hat\varphi(z(K)))^2\hat\varphi^{\prime\prime}(z(K))}>0,
\end{eqnarray*}
where the second equality follows by differentiating~\eqref{impl}.
Thus, $\hat P$ has a strictly positive second derivative on
$(\underline{K},\hat K)$, which completes the proof.~%
\end{pf}

\begin{remark*}
Note that $\hat{P}'(0+)
\geq0$ with equality if and only if
$\hat\varphi(0+)=\infty$.
\end{remark*}

We end this section by showing that $\hat\varphi$ can be recovered
directly from the put option prices $\hat P(K)$, at least on the domain
$(0,x_0]$. To do this, we define the function
$\varphi\dvtx(0,x_0]\to(0,\infty)$ by
%
\begin{equation}\label{invert}
\varphi(z)=\sup_{K:K\geq z}\frac{K-z}{\hat P(K)},
\end{equation}
where $\hat P$ is given by~\eqref{Px0}.

\begin{lemma}\label{dualitylemma}
\begin{longlist}[(b)]
\item[(a)] Suppose $f\dvtx(0,z_0]\rightarrow[1,\infty]$ is a nonnegative,
decreasing convex function on $(0,z_0]$ with $f(z_0)=1$ and
$f'(z_0)<0$. Define $g\dvtx(0,\infty) \rightarrow[0,\infty)$ by
%
\begin{equation}\label{g}
g(k) = \sup_{z:z\leq z_0} \frac{k-z}{f(z)}.
\end{equation}
\vspace*{-14pt}
\begin{enumerate}[(ii)]
\item[(i)]
$g(k)$ is then a nonnegative,
nondecreasing convex function with $(k-z_0)^+ \leq g(k) \leq k$ and
$g(k)=k-z_0$ for $k\geq k^*=z_0-1/f'(z_0)$.
\item[(ii)] $f$ and $g$
are self-dual in the sense that if, for $z\leq z_0,$ we define
\[
F(z) = \sup_{k:k\geq z} \frac{k-z}{g(k)},
\]
then $F \equiv f$ on $(0,z_0]$.
\end{enumerate}
\item[(b)] Similarly, assume that $g\dvtx(0,\infty)\to[0,\infty)$ is a
nonnegative, nondecreasing convex function with $(k-z_0)^+\leq g(k)\leq
k$ for all $k$. Also, assume that there exists a point $k^*>z_0$ such
that $g(k)=k-z_0$ for $k\geq k^*$ and $g(k)>k-z_0$ for $0\leq k<k^*$.
Define
\[
f(z)=\sup_{k:k\geq z}\frac{k-z}{g(k)}
\]
for $z\leq z_0$.
\begin{enumerate}[(ii)]
\item[(i)]
$f\dvtx(0,z_0]\to[0,\infty]$
is then a decreasing convex function with $f(z_0)=1$ and $f'(z_0)<0$.
\item[(ii)] $g$ and $f$ are self-dual in the sense that if we define
\[
G(k)=\sup_{z:z\leq z_0}\frac{k-z}{f(z)},
\]
then $G=g$ on $(0,\infty)$.
\end{enumerate}
\end{longlist}
\end{lemma}

\begin{pf}
See Appendix~\ref{app1}.\vadjust{\goodbreak}
\end{pf}

\begin{corollary}\label{transform}
 The function $\varphi$ coincides with the decreasing
fundamental solution $\hat\varphi$ on $(0,x_0]$.
\end{corollary}

\section{The inverse problem: The regular case}\label{reg}

We now consider the inverse problem. Let $P(K)$ be observed perpetual
put prices for all nonnegative values of the strike $K$. The idea is
that since $\hat\varphi$ satisfies the Black--Scholes equation
\eqref{ODE}, Corollary~\ref{transform} provides a way to recover the
volatility $\sigma(x)$ for $x\in(0,x_0]$ from perpetual put prices. In
this section we provide the details for the case where the observed put
prices are sufficiently regular. We assume that the observed put option
price $P\dvtx[0,\infty) \to[0,\infty)$ satisfies the following conditions
(cf. Propositions~\ref{simpleproperties} and~\ref{smoothproperties}
above).

\begin{hypothesis}\label{hyp}

\begin{enumerate}[(iii)]
\item[(i)]
$(K-x_0)^+\leq P(K)\leq K$ for all $K$.
\item[(ii)] There exists a
strike price $K^*$ such that $P(K)>(K-x_0)^+$ for all $K<K^*$ and
$P(K)=K-x_0$ for all $K\geq K^*$.
\item[(iii)] $P$ is continuously
differentiable on $(0,\infty)$ and twice continuously differentiable on
$(0,\infty)\setminus\{K^*\}$.
\item[(iv)] $P$ is strictly increasing on
$(0,\infty)$ with a strictly positive second derivative on $(0,K^*)$.
Moreover, $P''(K^*-):=\lim_{K\uparrow K^*}P''(K)$ exists and satisfies
$P''(K^*-)\in(0,\infty)$.
\end{enumerate}
\end{hypothesis}

Motivated by Corollary~\ref{transform}, we define the function
$\varphi\dvtx(0,x_0]\to(0,\infty)$ by
%
\begin{equation}\label{phi}
 \varphi(z)= \sup_{K:K\geq z}\frac{K-z}{P(K)}.
\end{equation}

\begin{proposition}\label{transform2}
The function $P$ can be recovered from $\varphi$ by
\[
P(K)=\sup_{z:z\leq z_0}\frac{K-z}{\varphi(z)}.
\]
\end{proposition}

\begin{pf}
This is a consequence of part (iii) of Lemma~\ref{dualitylemma}.
\end{pf}

\begin{proposition}\label{prop}
The function $\varphi\dvtx(0,x_0]\to(0,\infty)$ is twice
continuously differentiable with a positive second derivative, and it
satisfies $\varphi(x_0)=1$ and $\varphi^\prime(x_0)=-1/(K^*-x_0)$.
\end{proposition}

\begin{figure}

\includegraphics{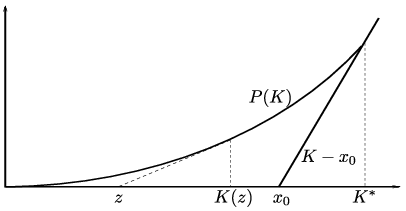}

\caption{For a given $z\leq x_0$ the value $\varphi(z)$ is given
by the slope of the tangent line to $P$ which passes through the point
$(z,0)$.}\label{fig2}
\end{figure}

\begin{pf}
For each $z\leq x_0$ there exists a unique $K=K(z)\in(z,K^*]$ for which
the supremum in~\eqref{phi} is attained. Geometrically, $K$ is the
unique value which minimizes the slope of the line through $(z,0)$ and
$(K,P(K))$ (cf. Figure~\ref{fig2}). Clearly, $K=K(z)$ satisfies the relation
%
\begin{equation}\label{Kz}
(K-z)P^\prime(K)=P(K).
\end{equation}
Reasoning as in the proof of Proposition~\ref{smoothproperties}, one
finds that $K(z)$ is continuously differentiable on $(0,x_0]$ with
%
\begin{equation}\label{palme}
\varphi^\prime(z)=\frac{-1}{P(K(z))}
\end{equation}
and $\varphi^\prime(x_0)=-1/(K^*-x_0)$. Differentiating~\eqref{palme}
with respect to $z$ gives
%
\begin{equation}\label{secder}
\varphi^{\prime\prime}(z)=
\frac{K^\prime(z)P^\prime(K(z))}{P^2(K(z))}=
\frac{(P^\prime(K(z)))^2}{(K(z)-z)P^2(K(z))P^{\prime\prime}(K(z))} ,
\end{equation}
where the second equality follows by differentiating~\eqref{Kz}. It
follows that $\varphi^{\prime\prime}(z)$ is continuous and positive,
which completes the proof.
\end{pf}
%

Next, we extend the function $\varphi$ to the whole positive real axis
so that $\varphi$ is convex, strictly positive, twice continuously
differentiable with a strictly positive second derivative and satisfies
$\varphi(\infty)=0$. We also define $\sigma^2(x)$ so that $\varphi$ is
a solution to the corresponding Black--Scholes equation, that is,
%
\begin{equation}\label{sigma}
\sigma^2(x)=2r\frac{\varphi(x)-x\varphi^\prime(x)}
{x^2\varphi^{\prime\prime}(x)}.
\end{equation}
Now, given this volatility function $\sigma(\cdot)$, we are in the
situation of Section~\ref{forward} and can thus define $\hat\varphi$ to
be the decreasing fundamental solution to the corresponding
Black--Scholes equation scaled so that $\hat\varphi(x_0)=1$. Moreover,
let $\hat P(K)$ be the corresponding perpetual put option price as
given by~\eqref{Px0}.

\begin{theorem}
Assume that Hypothesis~\ref{hyp} holds. The functions $\hat P$ and $P$
then coincide. Consequently, the volatility $\sigma(x)$ defined by
\eqref{sigma} solves the inverse problem.
\end{theorem}

\begin{pf}
Since the decreasing fundamental solution is unique up to a
multiplicative constant and $\varphi(x_0)=\hat\varphi(x_0)$, we have
$\varphi\equiv\hat\varphi$. Proposition~\ref{transform2} then yields
\[
\hat P(K)=\sup_{z:z\leq x_0}\frac{K-z}{\hat\varphi(z)}=
\sup_{z:z\leq x_0}\frac{K-z}{\varphi(z)}= P(K),
\]
which completes the
proof.
\end{pf}

\begin{remark*}
The inverse problem does not have
a unique solution. Indeed, there is
plenty of freedom when extending $\varphi$ (and thereby also $\sigma$)
for $x>x_0$. Note, however, that the volatility $\sigma$ is completely
determined by the given option prices for values below $x_0$.
\end{remark*}

We next show how to calculate the volatility that solves the inverse
problem directly from the observed option prices $P(K)$. To do that,
note that for each fixed $z\leq x_0$, the supremum in~\eqref{phi} is
attained at some $K=K(z)$ for which
%
\begin{equation}\label{z}
\varphi(z)=\frac{K-z}{P(K)},
\end{equation}
%
%
\begin{equation}
\varphi^\prime(z)=\frac{-1}{P(K)}
\end{equation}
and
%
\begin{equation}
\varphi^{\prime\prime}(z)=
\frac{(P^\prime(K))^2}{(K-z)P^2(K)P^{\prime\prime}(K)}
\end{equation}
[cf.~\eqref{palme} and~\eqref{secder}]. Since $\varphi$ satisfies the
Black--Scholes equation, we get
%
\begin{equation}\label{sig}
\sigma(z)^2z^2=2r\frac{\varphi(z)-z\varphi^\prime(z)}
{\varphi^{\prime\prime}} =\frac{2rK
P^2(K)P^{\prime\prime}(K)}{(P^\prime(K))^3}.
\end{equation}
Consequently, to solve the inverse problem we first determine $z$ by
\[
z=K-\frac{P(K)}{P^\prime(K)},
\]
and then, for this $z,$ we determines $\sigma(z)$ from~\eqref{sig}.

\section{The inverse problem: The irregular case}

Again, suppose we are given perpetual put prices $P(K)$ and a constant
interest rate $r>0$. Our goal is to construct a time-homogeneous
process which is consistent with the given prices. Unlike in the
regular case discussed in Section~\ref{reg}, we now impose no
regularity assumptions on the function $P$ beyond the necessary
conditions stated in Proposition~\ref{simpleproperties} and condition
(i) of Proposition~\ref{Kstar}. For a discussion of the necessity of
condition (i) of Proposition~\ref{Kstar}, see Section~\ref{noee}.

\begin{hypothesis}\label{hyp2}
\begin{enumerate}[(iii)]
\item[(i)]
For all $K$ we have $(K-x_0)^+\leq P(K)\leq K$.
\item[(ii)] $P$ is
nondecreasing and convex.
\item[(iii)] There exists $K^* \in
(x_0,\infty)$ such that $P(K)=K-x_0$ for $K\geq K^*$ and $P(K)>K-x_0$
for $K\in[x_0,K^*)$.
\end{enumerate}
\end{hypothesis}

\begin{theorem}\label{thm-irregular}
Given $P(K)$ satisfying Hypothesis~\ref{hyp2} and
given $r>0$, there exists a right-continuous (for $t>0$),
time-homogeneous Markov process $X_t$ with $X_0=x_0$ such that
\[
\sup_{\tau} \E^{x_0} [ e^{-r \tau} (K - X_\tau)^+ ] = P(K)\qquad
\forall K>0
\]
and such that $(e^{-rt}X_t)_{t \geq0}$ is
a local martingale.
\end{theorem}

\begin{remark*}
Although we wish to work in the
standard framework with
right-continuous processes, in some circumstances we have to allow for
an immediate jump. We do this by making the process right-continuous,
except possibly at $t=0$. At $t=0$ we allow a jump subject to the
martingale condition $\E[X_0]=x_0$.
\end{remark*}

Note that condition (iii) of Hypothesis~\ref{hyp2} excludes the
completely degenerate case where $P(x_0)=0$. If $P(x_0)=0,$ then,
necessarily, to preclude arbitrage, $P(K)= (K-x_0)^+$ and $X_t =
x_0e^{rt}$ is consistent with the prices $P(K)$. In this case
$\tau\equiv0$ is an optimal stopping time for every $K$.

Given $P(K)$ satisfying Hypothesis~\ref{hyp2}, we define $\varphi$ by
%
\begin{equation}\label{varphidef}
\varphi(x)= \sup_{K:K \geq x}\frac{K-x}{P(K)}
\end{equation}
for $x \in(0,x_0]$. For some values of $x$, the supremum in
\eqref{varphidef} may be infinite since $P$ may vanish on a nonempty
interval $(0,\underline K]$, where $\underline{K} = \sup\{K \dvtx \break P(K) = 0
\}$. By Lemma~\ref{dualitylemma}, $\varphi\dvtx(0,x_0]\to[1,\infty]$ is a
convex, decreasing, nonnegative function with $\varphi(x_0) = 1$.
Further,
%
\begin{equation}\label{leftderineq}
\varphi(x_0)-\varphi(x_0-\ep)\leq
1-\frac{K^*-x_0+\ep}{K^*-x_0} =\frac{-\ep}{K^*-x_0},
\end{equation}
so $\varphi'(x_0)\leq-1/(K^*-x_0)<0$. We define
\[
\underline{x} = \inf\{ x>0\dvtx \varphi(x)< \infty\},
\]
and in the case where $\underline{x}>0$ we see that $\varphi(x) =
\infty$ for $x < \underline{x}$. In fact, $\underline{x} >0$ if and
only if $\underline{K}>0$, and it is then easy to see that these two
quantities are equal.

We extend the definition of $\varphi$ to $(x_0,\infty)$ in any way
which is consistent with the convexity, monotonicity and nonnegativity
properties and such that $\lim_{x\to\infty}\varphi(x)=0$. It is
convenient to use $\varphi(x) = (x/x_0)^{\varphi'(x_0-)x_0}$, for
$\varphi'(x)$ is then continuous at $x_0$, and $\varphi$ is twice
continuously differentiable and positive on $(x_0,\infty)$.

Given $\varphi$, define $s\dvtx(\underline x,\infty) \mapsto
(-\infty,\infty)$ via
\[
s(x)= 2\int_{x_0}^x\varphi(y)\,dy +x_0-x\varphi(x)
\]
so that if $\varphi$ is differentiable we have $s'(x) = \varphi(x) - x
\varphi'(x)$. Then, $s$ is a concave, increasing function, which is
continuous on $(\underline{x},\infty)$. (It will turn out that $s$ is
the scale function, which explains the choice of label.) The function
$s$ has a well-defined inverse $g\dvtx(s(\underline
x),s(\infty))\to(\underline x,\infty)$, and if $s(\underline
x)>-\infty,$ then we extend the definition of $g$ so that
$g(y)=\underline x$ for $y\leq s(\underline x)$. Note that
$g\dvtx(-\infty,s(\infty)) \rightarrow[\underline{x},\infty)$ is a convex,
nondecreasing function with $g(0)=x_0$. Also, define $f(y) =
\varphi(g(y))$. Then, $f$ is decreasing and convex with
$f(0)=\varphi(x_0)=1$.

\begin{eg*}
For geometric Brownian motion we
have $s(0)=-x_0(1+\beta)/\break(1-\beta)$
and $s(\infty)=\infty$ if $\beta<1$, and $s(0)=-\infty$ and
$s(\infty)=x_0(1+\beta)/(\beta-1)$ if $\beta>1$. Moreover,
\begin{eqnarray*}
s(x) &=& x_0^\beta(x^{1 - \beta} - x_0^{1-\beta})(1+\beta)/(1-\beta),\\
g(y) &=& x_0\biggl[1 + \frac{y(1-\beta)}{x_0(1+\beta)}\biggr]^{1/(1-\beta)}
\end{eqnarray*}
and
\[
f(y) = \biggl[1 + \frac{y(1-\beta)}{x_0(1+\beta)}\biggr]^{-\beta/(1-\beta)}
\]
for $\beta\neq1$. If $\beta=1$, then the\vspace*{1pt} corresponding formulae are
$s(0)=-\infty$, $s(\infty)=\infty$, $s(x) = 2 x_0 \ln(x/x_0)$,
$g(y) =
x_0 e^{y/(2x_0)}$ and $f(y) = e^{-y/(2x_0)}$.
\end{eg*}

\begin{remark*}
Recall that a scale function is
only determined up to a linear
transformation. The choice $s(x_0)=0$ is arbitrary, but extremely
convenient, as it allows us to start the process $Z$, defined below, at
zero. The choice $s'(x) = \varphi(x) - x \varphi'(x)$ is simple, but a
case could be made for the alternative normalization $s'(x) =
(\varphi(x) - x \varphi'(x))/(1-x_0 \varphi'(x_0))$ for which
$s'(x_0)=1$. Multiplying $s$ by a constant has the effect of modifying
the construction defined in the next section, but only by the
introduction of a constant factor into the time changes. It is easy to
check that this leaves the final model $X_t$ unchanged.
\end{remark*}

Our goal is to construct a time-homogeneous process which is consistent
with observed put prices and such that $e^{-rt}X_t$ is a (local)
martingale. In the regular case we have seen how to\vadjust{\goodbreak} construct a
diffusion with these properties. We now have to allow for more general
processes, perhaps processes which jump over intervals, or perhaps
processes which have ``sticky'' points. One very powerful construction
method for time-homogeneous, martingale diffusions is via a time change
of Brownian motion, and it is this approach which we exploit.

\section{Constructing time-homogeneous processes as time changes
of Brownian motion} \label{sec-timechange}

In this section we extend the construction in Rogers and
Williams~\cite{RW2}, Section V.47, of martingale diffusions as time
changes of Brownian motion; see also It\^{o} and McKean~\cite{IM}, Section
5.1. The difference from the classical setting is that the
processes defined below may have ``sticky'' points and may jump over
intervals. Since diffusions are continuous by definition, the
resulting processes are not diffusions, but one might think of them as
``generalized diffusions'' (\cite{Kotani}, or ``gap diffusions''
\cite{Knight}), and they are ``as continuous as possible.''

Let $\nu$ be a Borel measure on $\R$ and let ${\mathbb
F}^B=({\mathcal
F}^B_u)_{u \geq0}$ be a filtration supporting a Brownian motion $B$
started at $0$ with local time process $L^z_u$. Define $\Gamma$ to be
the left-continuous increasing additive functional
%
\begin{equation}\label{Gammadef}
\Gamma_u = \int_{\R} L^z_u \nu(dz),\qquad\Gamma_0=0,
\end{equation}
and let $A$ be the right-continuous inverse of $\Gamma$, that is,
\[
A_t=\inf\{u\dvtx\Gamma_u>t\}.
\]
Note that $\Gamma$ is a nondecreasing process, so $A$ is well defined,
and $A_t$ is an $\mathbb F^B$-stopping time for each time $t$. Set
$Z_0=0$ and, for $t>0,$ set $Z_t = B_{A_t}$ and ${\mathcal
F}_t={\mathcal F}^B_{A_t}$. Note that $Z$ is right-continuous, except
possibly at $t=0$. The process $Z_t$ is a time-changed Brownian motion
adapted to the filtration ${\mathbb F}=({\mathcal F}_t)_{t \geq0}$ and
subject to mild nondegeneracy conditions on $\nu$ (see
Lemma~\ref{tc-lm} below), and the processes ${Z_t}$ and $Z_t^2 - A_t$
are local martingales. Further, if $\nu(dy) = dy/\gamma^2(y),$ then
$\Gamma_u = \int_0^u \gamma^{-2}(B_r)\,dr$ and $A_t = \int_{0}^{t}
\gamma^2(Z_s)\,ds$, so that $Z_t$ is a weak solution to $dZ_t =
\gamma(Z_t)\,dW_t$, and $Z$ is a diffusion in natural scale. Similarly,
if $\nu(dy)=dy/\gamma^2(y)$ in an interval, then $Z$ solves $dZ_t =
\gamma(Z_t)\,dW_t$ in this interval. The measure $\nu$ is called the
\textit{speed measure} of $Z$, although, as pointed out by Rogers and
Williams, $\nu$ is large when $Z$ moves slowly.

The measure $\nu$ may have atoms, and it may have intervals on which it
places no mass. If there is an atom at $\hat{z,}$ then $d
\Gamma_u/du>0$ whenever $B_u=\hat{z}$, and then the time-changed
process is ``sticky'' there. Conversely, if $\nu$ places no mass in
$(\alpha,\beta),$ then $\Gamma$ is constant on any time-periods that
$B$ spends in this interval, and the inverse time change $A$ has a
jump. In particular, $Z_t$ spends no time in this interval. If
$\nu(\{\tilde{z}\})=\infty,$ then $\Gamma_u=\infty$ for any $u$ greater
than the first hitting time $H^B_{\tilde{z}}$ by $B$ of level
$\tilde{z}$. In that case, $A_\infty\leq H^B_{\tilde{z}}$ so that if
$Z$ hits $\tilde{z}$, then $\tilde{z}$ is absorbing for $Z$. The other
possibility is that $Z$ tends to this level without reaching it in
finite time.\vadjust{\goodbreak}

Define $\overline{z}_\nu\in(0,\infty]$ and $\underline{z}_\nu\in
[-\infty,0)$ via
\begin{eqnarray*}
\overline{z}_\nu&=& \inf\{ z>0\dvtx \nu((0,z])=\infty\}\quad\mbox{and}\\
\underline{z}_\nu&=& \sup\{ z<0\dvtx \nu([z,0))=\infty\}.
\end{eqnarray*}
The cases where $\overline{z}_\nu=0$ or $\underline{z}_\nu=0$
correspond to the degenerate case $X_t = x_0e^{rt}$ mentioned in the
previous section, and we exclude them. The following lemma provides a
guide to sufficient conditions for a time change of Brownian motion to
be a local martingale and therefore provides insight into the
constructions of local martingales via time change that we develop in
the next section.

\begin{lemma}\label{tc-lm}
Suppose that either $\overline{z}_\nu<\infty$ or $\nu$ charges
$(a,\infty)$ for each $a$, and, further, that either
$\underline{z}_\nu>-\infty$ or $\nu$ charges $(-\infty,a)$ for each $a$.
Then, $Z_t=B_{A_t}$ is a local martingale.
\end{lemma}
\begin{pf} See Appendix~\ref{app2}.
\end{pf}

\section{Constructing the model}
We now show how to choose the measure $\nu$ which gives the process we
want. Define $\nu$ via
%
\begin{equation}\label{nudef}
 \nu(dy) = \frac{1}{2r} \frac{g''(dy)}{g(y)},
\end{equation}
where $g''(dy)$ is the measure defined by the second order distribution
derivative of $g$, and let $\nu(\{y\}) = \infty$ for $y\leq s(0)$ in
the case where $\underline x=0$ and $s(0)> -\infty$. Similarly, in the
case where $s(\infty)<\infty,$ we set $\nu(\{y\})=\infty$ for
$y\geq
s(\infty)$. Where $g'$ is absolutely continuous it follows that $\nu$
has a density with respect to Lebesgue measure, but, more generally,
\eqref{nudef} can be interpreted in a distributional sense.

Now, for this $\nu$ we can use the construction of the previous section
to give a process $Z_t$. If we set $X_t = g(Z_t)$, then, subject to
the hypotheses of Lemma~\ref{tc-lm}, $Z_t=s(X_t)$ is a local
martingale, so that $s$ is a scale function for $X$. The process $X$ is
our candidate process for which the associated put prices are given by
$P$.

\begin{eg*}
For geometric Brownian motion,
\[
\frac{\nu(dy)}{dy} = \frac{1}{2r} \frac{g''(y)}{g(y)} =
\frac{1}{2rx_0^2} \frac{\beta}{(1+\beta)^2} \biggl[ 1 + \frac{y
(1-\beta)}{x_0(1+\beta)} \biggr]^{-2}
\]
for $y\in(s(0),s(\infty))$.
In the case $\beta=1$ this simplifies to
\[
\frac{\nu(dy)}{dy} = \frac{1}{8 x_0^2 r} =
\frac{1}{C^2},
\]
where $C= 2 x_0 \sqrt{2r}$. Then, $\Gamma_u = u
C^{-2}$, $A_t = t C^2$, $Z_t = B_{t C^2} \stackrel{\mathcal D}{=} C
\tilde{B}_t$ for a Brownian motion $\tilde{B}$ and $X_t =
x_0e^{Z_t/2x_0} = x_0 e^{\sqrt{2r}\tilde{B}_t}$.
\end{eg*}

\begin{remark*}
If the put price $P(K)$ satisfies
the regularity conditions of
Hypothesis~\ref{hyp}, then the scale function $s$ and its inverse $g$
are $C^2$ and satisfy
\[
g'(y)s'(g(y))=1
\]
and
\[
g''(y)s'(g(y))+(g'(y))^2s''(g(y))=0.
\]
Moreover,
$\sigma^2(x)x^2 s''(x) + 2rxs'(x)=0$ so that
%
\begin{equation}\label{geqn}
\frac{g(y)^2 \sigma(g(y))^2 g''(y)}{2 g'(y)^2} = r g(y) .
\end{equation}

Consequently, the speed measure $\nu$ is given by
\begin{eqnarray*}
\nu(dy) &=& \frac{1}{2r}\frac{g''(y)}{g(y)}\,dy\\
&=& \frac{1}{2r}\frac{-(g'(y))^2s''(g(y))}{s'(g(y))g(y)}\,dy\\
&=& \frac{(g'(y))^2}{\sigma^2(g(y))g^2(y)}\,dy,
\end{eqnarray*}
and the diffusion $Z$ is the solution to the stochastic differential
equation
\[
dZ_t=\frac{\sigma(g(Z_t))g(Z_t)}{g'(Z_t)}\,dW_t.
\]
Applying
It\^o's formula to $X_t=g(Z_t)$ yields
\begin{eqnarray*}
dX_t &=& g'(Z_t)\,dZ_t+ \frac{1}{2}g''(Z_t)(dZ_t)^2\\
&=& \frac{\sigma^2(X_t)X_t^2 g''(Z_t)}{2(g'(Z_t))^2}\,dt+
\sigma(X_t)X_t\,dW_t\\
&=& rX_t\,dt + \sigma(X_t)X_t\,dW_t,
\end{eqnarray*}
where we use~\eqref{geqn} for the final equality. We thus recover the
diffusion model from the regular case described in Section~\ref{reg}.
\end{remark*}

Recall that $\Gamma_u = \int_{\R} L^z_u \nu(dz)$ and let $\xi$ be the
first explosion time of $\Gamma$. Note that by construction $\Gamma$ is
continuous for $t < \xi$ and left-continuous at $t = \xi$. Since~$\nu$
is infinite outside the interval $[s(0),s(\infty)],$ we also have the
expression $\xi= \inf\{ u \dvtx B_u \notin[s(0),s(\infty)] \} =
H^B_{s(0)} \wedge H^B_{s(\infty)}$. The inverse scale function $g$ is
convex on $(s(0), s(\infty))$, but may have a jump (from a finite to an
infinite value) at~$s(\infty)$. In that case we take it to be
left-continuous at $s(\infty)$ so that we may have $\bar{g}:= \lim_{z
\uparrow\infty}g(s(z))$ is finite.

For $0 \leq u < \xi$, define $M_u = e^{-r \Gamma_{u}} g(B_{u})$ and
$N_u = e^{-r \Gamma_{u}} f(B_{u})$.

\begin{lemma}\label{lem-lm}
$M=(M_u)_{0 \leq u < \xi}$ and $N=(N_u)_{0 \leq u <
\xi}$ are ${\mathbb F}^B$-local martingales.
\end{lemma}

\begin{pf*}{Sketch of proof} Suppose that $\varphi$ is twice
continuously differentiable with a positive second derivative. Then,
$g$ is twice continuously differentiable. For $u<\xi$, applying
It\^{o}'s formula to $M_u = e^{-r \Gamma_u} g(B_u)$ gives
\[
e^{r \Gamma_u}\,dM_u = g'(B_u)\,dB_u + \biggl[ -r \frac{d \Gamma_u}{du}
g(B_u) + \frac{1}{2} g''(B_u) \biggr]\, du.
\]
But, by definition,
$d\Gamma_u/du = g''(B_u)/(2r g(B_u))$, so $M$ is a local martingale, as
required.

A similar argument can be provided for the process $N$. For the general
case, see the \hyperref[app]{Appendix}.
\end{pf*}

Since $M$ and $N$ are nonnegative local martingales on $[0,\xi),$ they
converge almost surely to finite values, which we label $M_{\xi}$ and
$N_\xi$.
In particular,\vspace*{-2pt} if $\xi= H^B_{s(0)}$, then $M_{\xi}=0$. However, if
$\xi= H^B_{s(\infty)},$ then there are several cases. The fact that a
nonnegative local martingale converges means that we cannot have both
$\Gamma_\xi<\infty$ and $\bar{g}=\lim_{z \uparrow\infty} g(s(z)) =
\infty$. Instead, if $\Gamma_\xi<\infty,$ then $\bar{g}<\infty$ and
$M_{\xi}=e^{-r \Gamma_\xi} \bar{g}$. If $\Gamma_\xi= \infty$ and
$\bar{g}<\infty$, then $M_{\xi}=0$, whereas if $\Gamma_\xi= \infty$
and $\bar{g}=\infty$, then $(M_u)_{u < \xi}$ typically has a nontrivial
limit. Similar considerations apply to $N$.

Recall that $A$ is the right-continuous inverse of $\Gamma$ and define
the time-changed processes $\tilde{M}_t = M_{A_t}$ and $\tilde{N}_t =
N_{A_t}$. Note that these processes are adapted to ${\mathbb F}$ and
that, at least for $t < \Gamma_\xi$, we have $\Gamma_{A_t}=t$,
$\tilde{M}_t = e^{-rt} g(Z_t) = e^{-rt}X_t$ and $\tilde{N}_t =
e^{-rt}f(Z_t) = e^{-rt}\varphi(X_t)$.

If $(s(0),s(\infty))=\R,$ then $\xi= \infty$, $\Gamma_\xi=\infty
$ and
$\tilde{M}$ is defined for all $t$.

If $s(0)>-\infty,$ then we may have $\xi= H^B_{s(0)}$. In this\vspace*{-2pt} case
either $\Gamma_\xi= \Gamma_{H^B_{s(0)}}=\infty$, whence $\tilde
{M}$ is
defined for\vspace*{-1pt} all $t$ as before, or $\Gamma_{\xi}<\infty$. Then,
$\tilde{M}_{\Gamma_\xi} = M_\xi= e^{-r \Gamma_\xi} g(B_{\xi}) = 0,$
and we set $\tilde{M}_t = 0$ for all $t>\Gamma_\xi$. It follows that
$X_t=0$ for all $t\geq\Gamma_\xi$, and 0 is an absorbing state.

Similarly, if $s(\infty)<\infty,$ then\vspace*{-1pt} we may have $\xi=
H^B_{s(\infty)}$. Then, either $\Gamma_{\xi}=\infty$, whence
$\tilde{M}$ is defined\vspace*{-1pt} for all $t$, or $\Gamma_{\xi}<\infty$. In the
latter case, if $\xi= H^B_{s(\infty)}<\infty$ and
$\Gamma_{\xi}<\infty,$ then $\tilde{M}_{\Gamma_{\xi}} = M_{\xi}
= e^{-r
\Gamma_\xi} \bar{g}$. We set\vspace*{-1pt} $\tilde{M}_t = \tilde{M}_{\Gamma_\xi
}$ for
all $t>\Gamma_\xi$, and it follows that for $t>\Gamma_\xi$, $X_t :=
e^{rt} \tilde{M}_t = e^{r(t - \Gamma_\xi)} \bar{g}$. Thus, for $t >
\Gamma_\xi$, $X$ grows deterministically. An example of this situation
is given in Example~\ref{eg-gbar<infty} below. (In fact, the case where
$\bar{g}<\infty$, which depends on the behavior of the scale function~$s$ to the right of $x_0$, can always be avoided by a suitable choice
of the extension to $\varphi$.)

We want to show how $\tilde{M}$ and $X$ inherit properties from $M$.
The key idea below is that, loosely speaking, a~time change of a
martingale is again a martingale. Of course, to make this statement
precise we need strong control on the time change. (Without such\vadjust{\goodbreak}
control the resulting process can have arbitrary drift. Indeed, as
Monroe~\cite{Mon} has shown, any semimartingale can be constructed from
Brownian motion via a time change.) We have the following result, the
proof of which is given in the \hyperref[app]{Appendix}.

\begin{corollary}\label{e-rtX=lm}
The process $(e^{-rt}X_t)_{t \geq0}$ is a local martingale.
\end{corollary}

We can perform a similar analysis on $N$ and $\tilde N$ and use similar
ideas to ensure that $\tilde N$ is defined on $\mathbb{R}_+$. The proof
that $\tilde{N}$ is a local martingale mirrors that of
Corollary~\ref{e-rtX=lm}.

\begin{corollary}\label{e-rtf=lm}
The process $(e^{-rt}\varphi(X_t))_{t \geq0}$ is a local martingale.
\end{corollary}

\section{Determining the put prices for the candidate process}
Recall the definitions of $s$, $g$ and $\nu$ via $s'(x) = \varphi(x) -
x \varphi'(x)$, $g \equiv s^{-1}$ and $\nu(dy) = g''(dy)/(2rg(y))$.
Suppose that $Z$ is constructed from $\nu$ and a Brownian motion using
the time change $\Gamma$ and construct the candidate price process via
$X_t = g(Z_t)$. By Corollary~\ref{e-rtX=lm}, the discounted price
$e^{-rt}X_t$ is a (local) martingale. To complete the proof of
Theorem~\ref{thm-irregular} we need to show that for the candidate
process $X_t$, the function
\[
\hat{P}(K) := \sup_{\tau} \E[e^{-r \tau}(K- X_\tau)^+]
\]
is such that $\hat{P}(K) \equiv P(K)$ for all $K \geq0$.

Unlike the regular case, the process $X$ that we have constructed may
have jumps. For this reason, for $x<x_0$ we modify the definition of
the first hitting time so that $H_x = \inf\{ u>0\dvtx X_u \leq x \}$.

\begin{theorem}
The perpetual put prices for $X$ are given by $P$.
\end{theorem}

\begin{pf}
Fix $x\in(\underline{x},x_0)$. Suppose first that
$x$ is such that $\Gamma$ is strictly increasing whenever the Brownian
motion $B$ takes the value $s(x)$. Then, $X_{H_x}=x$. More generally,
the same is true whenever $\nu((s(x)-\delta,s(x)])>0$ for every
$\delta>0$. By Corollary~\ref{e-rtf=lm} we have that $(e^{-rt}
\varphi(X_t))_{t \leq H_x}$ is a local martingale,\vspace*{1pt} and $\varphi$ is
bounded on $[x,\infty),$ so it follows that $e^{-r(t \wedge
H_x)}\varphi(X_{t \wedge H_x})$ is a bounded martingale and
$\varphi(x_0) = \E^{x_0}[ e^{-rH_x} \varphi(x)]$. Hence,
\[
\hat{P}(K) \geq\E^{x_0}[ e^{-rH_x} (K-x)] = (K-x)
\frac{\varphi(x_0)}{\varphi(x)}=\frac{K-x}{\varphi(x)}.
\]

Otherwise, fix ${x}^-(x) = \inf\{ w<x\dvtx \nu((s(w),s(x)]) = 0 \}$ and
$x^+(x)=\break\sup\{ w>x\dvtx \nu([s(x),s(w))=0 \}$. It must be the case that
$\varphi$ is linear on $(x^{-}(x),x^+(x))$, bounded on
$[x^{-}(x),\infty)$ and
\begin{eqnarray*} \hat{P}(K) &\geq& \max_{w \in\{x^-,x^+\}} \E^{x_0}[
e^{-rH_w} (K-w)] \\
& = & \max_{w \in\{x^-,x^+\}} \frac{K-w}{\varphi(w)} \geq
\frac{K-x}{\varphi(x)} .
\end{eqnarray*}
It follows that
%
\begin{equation}\label{waldner}
\hat{P}(K) \geq\sup_{x:x\leq x_0}
\frac{K-x}{\varphi(x)} = P(K).
\end{equation}
[Clearly, if $x<\underline{x}$, then $(K-x)/\varphi(x)=0$, so the
supremum cannot be attained for such an $x$.]

To prove the reverse inequality, we first claim that the left
derivative $D^-\varphi$ of the convex function $\varphi$ satisfies
%
\begin{equation}\label{leftder}
D^-\varphi(x_0):=\lim_{\ep\downarrow
0}\frac{\varphi(x_0)-\varphi(x_0-\ep)}{\ep}=\frac{-1}{K^*-x_0}.
\end{equation}
To prove~\eqref{leftder}, first note that it follows from
\eqref{leftderineq} that
\[
D^-\varphi(x_0)\leq\frac{-1}{K^*-x_0}.
\]
Conversely, note that for each $\delta>0$ there exists a nonempty
interval $(x_0-\ep,x_0),$ on which
%
\begin{equation}\label{leftderineq2}
\varphi(x)\leq
\frac{K^*-\delta-x}{K^*-\delta-x_0}.
\end{equation}
To see this, let $\delta>0$ be small and draw the tangent line to $P$
that passes through the point $(K^*-\delta,K^*-x_0-\delta)$. Let
$x_0-\ep$ be the $x$-coordinate of the point of intersection between
the tangent line and the $x$-axis. Then, for all $x\in[x_0-\ep,x_0]$ we
have that the line through $(x,0)$ and $(K^*-\delta,K^*-x_0-\delta)$ is
below the graph of $P$. Consequently,~\eqref{leftderineq2} holds.
Therefore, for $x\in(x_0-\ep,x_0)$ we have
\[
\frac{\varphi(x_0)-\varphi(x)}{x_0-x}\geq
\frac{-1}{K^*-\delta-x_0}.
\]
Thus, $D^-\varphi(x_0)\geq-1/(K^*-x_0)$
since $\delta>0$ is arbitrary, so~\eqref{leftder} follows.

We next claim that for each fixed $K\leq K^*$ we have
%
\begin{equation}\label{obama}
\varphi(x) \geq(K-x)^+/P(K)
\end{equation}
for all $x$. Clearly, this holds for $x\geq K$ and for $x\leq x_0$.
Similarly, if $x_0 < x< K$, then it follows from~\eqref{leftder} and
the convexity of $\varphi$ that
\[
\varphi(x)\geq\frac{K^*-x}{K^*-x_0}\geq\frac{K-x}{K-x_0}
\geq\frac{K-x}{P(K)}.
\]
It follows from~\eqref{obama} and
Corollary~\ref{e-rtf=lm} that for any stopping rule $\tau$ we have
\[
\E^{x_0}[e^{-r \tau} (K-X_\tau)^+ ] \leq P(K)\E^{x_0}[e^{-r \tau}
\varphi(X_\tau) ] \leq P(K)\varphi(x_0) = P(K).
\]
Hence, $\hat P(K)
\leq P(K)$ for $K\leq K^*$ and, in view of~\eqref{waldner}, $\hat
P(K)=P(K)$.

For $K> K^*$ it follows from $\hat P(K^*)=P(K^*)=K^*-x_0$, the
convexity of $\hat P$ and Hypothesis~\ref{hyp2} that $\hat
P(K)=K-x_0=P(K)$, which completes the proof.
\end{pf}



\section{Examples}\label{sectionexamples}
The following examples illustrate the
construction of the previous sections. The list of examples is not
intended to be exhaustive, but rather indicative of the types of
behavior that can arise. In each example we assume $x_0=1$.

\subsection{The smooth case}
We have studied the case of exponential Brownian motion throughout. It
is very easy to generate other examples, for example, by choosing a
smooth decreasing convex function [with $\varphi(x_0)=1$ and $\lim_{x
\uparrow\infty}\varphi(x)=0$] and defining other quantities from
$\varphi$.
\begin{eg}
Suppose $\varphi(x) = (x+1)/(2x^2)$. Then, from~\eqref{ODE} we obtain
\[
\sigma^2(x) = r \frac{2x + 3}{x + 3},\qquad x>0,
\]
and from~\eqref{Px0},
\[
P(K) = \frac{(K+9)^{3/2}(K+1)^{1/2} - (27 +18K-K^2)}{4},\qquad
K \leq5/3,
\]
with $P(K)=(K-1)$ for $K \geq5/3$.
\end{eg}

\subsection{Kinks in $P$}\label{kink}
If the first derivative of $P$ is not continuous, then we
find that~$\varphi$ is linear over an interval $(\alpha,\beta)$, say.
Then, $s'$ is constant on this interval and $g$ is linear over the
interval $(s(\alpha),s(\beta))$. It follows that $\nu$ does not charge
this interval, so $\Gamma_u$ is constant whenever $B_u \in
(s(\alpha),s(\beta))$, and $A_t$ has a jump. $Z_t$ then jumps over the
interval $(s(\alpha),s(\beta))$, and $X_t$ spends no time in
$(\alpha,\beta)$.

\begin{eg}\label{eg:kink}
Suppose that $P(K)$ satisfying
Hypothesis~\ref{hyp2} is given by
\[
P(K) = \cases{
{K^2}/{8,}
&\quad $0 < K \leq27/32$, \cr
4K^3/27, &\quad $27/32 \leq K \leq3/2$, \cr
(K-1), &\quad $3/2 \leq K$.
}
\]
$P$ is then continuous, but $P'$ has a jump at $K=27/32$.

Using~\eqref{varphidef} we find that
\[
\varphi(x) = \cases{
{2}x^{-1}, &\quad $0 < x \leq27/64$, \cr
x^{-2}, &\quad $x>9/16$
}
\]
(strictly speaking, there is some freedom in the choice of $\varphi$
for $x \geq x_0 \equiv1$, but the power function $x^{-2}$ is a natural
choice). Over the region $I=[27/64,9/16],$ $\varphi$ is given by linear
interpolation. The corresponding scale function is linear on $I$ and in
the construction of $Z$, $\nu$ assigns no mass to $s(I)$. The\vspace*{1pt} process
$X$ is a generalized diffusion with diffusion coefficient given by
$\sigma(x)= \sqrt{2r}$ for $x \leq27/64$, $\sigma(x)=\sqrt{r}$ for $x
\geq9/16$ and $\sigma(x) = \infty$ for $x \in I$.
\end{eg}

\subsection{Linear parts to $P$}
In this case, the derivative of $\varphi(x)$ is discontinuous at a
point $\gamma,$ say. Then, $s'$ is also discontinuous at this point,
and $g'$ is discontinuous at $s(\gamma)$. It follows that $\nu$ has a
point mass at $s(\gamma)$, and that $\Gamma_u$ includes a multiple of
the local time at $s(\gamma)$.

\begin{eg}
Suppose that$P(K)$ satisfying Hypothesis~\ref{hyp2} is given
by
\[
P(K) = \cases{
{8K^3}/{27,}
&\quad $0 < K \leq3/4$, \cr
(2K-1)/4, &\quad $3/4 \leq K \leq1$, \cr
K^2/4, &\quad $1 \leq K \leq2$, \cr
(K-1), &\quad $2 \leq K$.
}
\]
$P$ is then convex, but is linear on the interval $[3/4,1]$.
We have
\[
\varphi(x) = \cases{
x^{-2}/2,
&\quad 0 < x $\leq1/2$, \cr
x^{-1}, &\quad $x>1/2$,
}
\]
where we have chosen to extend the definition of $\varphi$ to
$(1,\infty)$ in the natural way. Then, $s(x) = 3 - 2 \ln2 - 3/2x$ for
$x<1/2$ and $s(x)=2 \ln x$ otherwise. It follows that $g$ is everywhere
convex, but has a discontinuous first derivative at $z=-2 \ln2$, and
that the corresponding measure $\nu$ has a positive density with
respect to Lebesgue measure \textit{and} an atom of size $r^{-1}/12$ at
$-2 \ln2$. In the terminology of stochastic processes, the process $Z$
is ``sticky'' at this point; for a~discussion of sticky Brownian
motion, see Amit~\cite{Ami} or, for the one-sided case, see
Warren~\cite{War}.
\end{eg}

If $P$ is piecewise linear (e.g., if $P$ is obtained by linear
interpolation from a finite number of options), then $\varphi$ is
piecewise linear, $s$ is piecewise linear, $g$ is piecewise linear and
$\nu$ consists of a series of atoms. As a consequence the process $Z_t$
is a continuous-time Markov process on a countable state space [at
least while $Z_t < s(x_0)\equiv0$], in which transitions are to
nearest neighbors only. Holding times in states are exponential and the
jump probabilities are such that $Z_t$ is a martingale.

In turn this means that $X_t$ is a continuous-time Markov process on a
countable set of points (at least while $X_t < x_0$).

\begin{eg}\label{eg-gbar<infty}
Suppose
\[
P(K)=\cases{
K/3, &\quad $K \leq1$,\cr
(2K-1)/3, &\quad $1 \leq K \leq2$,\cr
(K-1), &\quad $K \geq2$.
}
\]
This is consistent with a situation in which only two perpetual
American put options trade, with strikes $1$ and $3/2,$ and prices
$1/3$
and $2/3$, in which case we may assume that we have extrapolated from the
traded prices to a put pricing\vadjust{\goodbreak} function~$P(K)$ which is consistent with
the traded prices. The function
\[
\varphi(x)=\cases{
3-3x, &\quad $x<1/2$,\cr
2-x, &\quad $1/2\leq x\leq2$,\cr
0, &\quad $x> 2$
}
\]
is a possible choice of $\varphi$. Then,
\[
s(x)=\cases{
3x-5/2, &\quad $x<1/2$,\cr
2x-2, &\quad $1/2\leq x\leq2$,\cr
2, &\quad $x> 2$.
}
\]
The inverse of $s$ is given by
\[
g(y)=\cases{
y/3+5/6, &\quad $-5/2\leq y<-1$,\cr
y/2+1, &\quad $-1\leq y\leq2$.
}
\]
The corresponding measure $\nu$ assigns no mass to the intervals
$(-5/2,-1)$ and $(-1,2),$ but has a point mass of size $1/(6r)$ at
$-1$. The corresponding process $X$ has state space $\{0\} \cup\{1/2\}
\cup
[2,\infty)$ and is such that:
\begin{itemize}
\item at $t=0+$, $X$ jumps to $1/2$ or 2 with probabilities
$2/3$ and $1/3$, respectively;
\item if ever $X_{t_0}\geq2,$ then $X_t = X_{t_0}
e^{r(t-t_0)}$ thereafter;
\item zero is an absorbing state for $X$.
\end{itemize}
To examine what happens if $X$ ever reaches $1/2$, note that $\xi=
H_{-5/2} \wedge H_2$ and $\Gamma_{\xi} = (1/6r) L^{-1}_{\xi}$ (where
the superscript denotes local time at $-1$ rather than an inverse) and
then
\[
\Prob(A_t<\xi) = \Prob\biggl(t< \frac{1}{6r}L^{-1}_{\xi}\biggr)
=\int_{6rt}^\infty\frac{1}{2}e^{-y/2}\,dy=e^{-3rt},
\]
where we have
used the known density of $L^{-1}_{H_{-5/2} \wedge H_2}$ (cf. page 213
in~\cite{BS}). This implies that if $X$ ever reaches $1/2$, then it
stays there for an exponential length of time, rate $3r$, and jumps to
2 with probability $1/3$ and zero with probability $2/3$.

Note that for the continuous-time Markov process $X_t$, conditional on
$X_t=1/2,$ we have
\[
\lim_{\Delta\downarrow0} \frac{1}{\Delta} \E[X_{t+\Delta}-X_t] =
3r \biggl[ \frac{1}{3} (2 - X_t) + \frac{2}{3} (-X_t) \biggr]
= \frac{r}{2} = rX_t.
\]
Also, for this process
\begin{eqnarray*}
P(K) & = & \max_{\tau= 0, H_{1/2}, H_0 } \E[ e^{- r \tau} (K-
X_\tau)^+]
\\
& = & \max\{ K-1, (2K-1)/3, K/3 \},
\end{eqnarray*}
so we recover the put price function given at the start of the example.
\end{eg}


\subsection{Positive gradient of $P$ at zero [i.e.,  $P'(0)>0$]}
In this case\break $\lim_{x \downarrow0} \varphi(x)<\infty$. It follows that
$s(0)>-\infty$ and the resulting diffusion $X_t$ can hit zero in finite
time. Recall that the diffusion $X$ is constructed so that 0 is an
absorbing endpoint.

\begin{eg}
Suppose that
\[
P(K)=\cases{
K/2, &\quad $K<1$, \cr
(K+1)^2/8, &\quad $1 \leq K \leq3$,\cr
K-1, &\quad $K \geq3$.
}
\]
Then, $\varphi(x)=2(x+1)^{-1}$
and $x^2\sigma(x)^2 = r(x+1)(2x+1)$ so that $dX_t = rX_t\,dt +
\sqrt{r(X_t+1)(2X_t+1)}\,dB_t$.
\end{eg}

The following example covers the case of mixed linear and smooth parts
of~$P(K)$ and shows an example where reflection, local times and jumps
all form part of the construction.

\begin{eg}
Suppose $P(K)$ satisfies
\[
P= \cases{
K/4, &\quad $K\leq1$,\cr
K^2/4, &\quad $1 \leq K \leq2$,\cr
(K-1), &\quad $K\geq2$.
}
\]
Note that $P'$ has a jump at
$K=1$. We have $\varphi(x) =4 - 4x$ for $x < 1/2$ and $\varphi(x)=1/x$
for $1/2\leq x \leq1$. We assume this formula also applies on
$[1,\infty)$. Then,
\[
s(x)=\cases{
2\ln x, &\quad $x\geq1/2$,\cr
4x-2-2\ln2, &\quad $x< 1/2$
}
\]
and
\[
g(y)=\cases{
e^{y/2}, &\quad $y\geq-2\ln2$,\cr
y/4 + (1+\ln2)/2, &\quad $-2-2\ln2< y< -2\ln2$.
}
\]
Consequently, $\nu(dy) = \frac{1}{8r}\,dy$ for $y\geq-2\ln2$, no mass
is assigned to the interval $(-2-2\ln2, -2\ln2)$ and $v(\{
y\})=\infty$ for $y\leq-2-2\ln2$. It follows that for $x \geq1/2$ we
have $\sigma^2(x) = 2r$, and then
%
\begin{equation}\label{8.5bexample}
dX_t = rX_t\,dt + \sqrt{2r} X_t\,dB_t,
\end{equation}
at least until the first hitting time of $1/2$. To allow for behavior
at $1/2$ the general construction includes a local time reflection
\textit{and} a compensating downward jump are added at instants when $X_t=1/2$.
The jump takes the process to zero, where it is absorbed.

Alternatively, the process can be formalized as follows. Let $I_t$ be
the infimum process given by $I_t= - \inf_{u \leq t} \{ (B_u + \ln
2/\sqrt{2r}) \wedge0 \}$. By Skorokhod's lemma, $B_t + I_t$ is then a
reflected Brownian motion [reflected at the level $- (\ln
2/\sqrt{2r})$] and $e^{\sqrt{2r} (B_t + I_t)} \geq1/2$.

Let $N^{\lambda}$ be a Poisson process with rate $\lambda$, independent
of $B$, and let $T^{\lambda}$ be the first event time. The compensated
Poisson process $(N^{\lambda}_t - \lambda t)_{t \geq0}$ and the
compensated Poisson process stopped at the first jump $(N^{\lambda}_{t
\wedge T^{\lambda}} - \lambda( t \wedge T^{\lambda}))_{t \geq0}$ are
then martingales. The time change $(N^{\lambda}_{I_t \wedge
T^{\lambda}} - \lambda( I_t \wedge T^{\lambda}))_{t \geq0}$ is also a
martingale.

Take $\lambda= \sqrt{2r}$ and define $X$ via $X_0=1$ and
\[
dX_t = rX_t\,dt + \sqrt{2r} X_t \biggl( dB_t + dI_t -
\frac{dN^{\sqrt{2r}}_{I_t}}{\sqrt{2r}} \biggr),\qquad t \dvtx I_t
\leq T^{\sqrt{2r}}.
\]
Note that at the first jump time of the
time-changed Poisson process, $X$ jumps from~$1/2$ to zero.

By construction, $(e^{-rt}X_t)_{t \geq0}$ is a martingale.
\end{eg}

\subsection{$P$ is zero on an interval}
Now, consider the case where $P(K)=0$ for $K \leq\underline{K}$. We
then find that $\varphi(x) = \infty$ for $x \leq\underline{x}$, where
$\underline{x} = \underline{K}$. Depending on whether the right
derivative $P'(\underline{K}+)$ is zero or positive,
$\varphi(\underline{x}+)$ may be infinite or finite. In the former case
we have that $X_t$ does not reach $\underline{x}$ in finite time. In
the latter case $X_t$ does hit $\underline{x}$ in finite time.

The first example is typical of the case where
$\varphi(\underline{x}+)=\infty$ or, equivalently, where there is
smooth fit of $P$ at $\underline{K}$.

\begin{eg}
Suppose $X_0=1$ and that $P(K)$ solves
\[
P(K) = \cases{
0, &\quad $K \leq1/2$,\cr
(2K-1)^2/8, &\quad $1/2\leq K\leq3/2$,\cr
K-1, &\quad $K \geq3/2$.
}
\]
$P'$ is then continuous and for $1/2<x<1$ we have $\varphi(x) = (2x-1)^{-1}$.
We also have $\varphi(x) = \infty$ for $x \leq1/2$. As usual, there is
some freedom when extending $\varphi$ to $(1,\infty)$, but for
definiteness we assume that the formula $\varphi(x) = (2x-1)^{-1}$
applies there as well.

It follows that $\eta(x)^2 \equiv(x \sigma(x))^2 = r(2x-1)(4x-1)/4$.
Note that since $\varphi(1/2) = \infty$ we have that $H_{1/2}$ (the
first hitting time of $1/2$) is infinite. Hence,
\[
dX_t = rX_t\,dt + \sqrt{\frac{r(2X_t-1)(4X_t-1)}{4}}\,dB_t,
\qquad t \leq H_{1/2},
\]
is consistent with the observed put
prices, and since the process never hits $1/2$ it is not necessary to
describe it beyond $H_{1/2}$.
\end{eg}

Now, consider the other case where $P'(\underline{K})>0$.
\begin{eg}
Suppose $X_0=1$ and that $P(K)$ solves
\[
P(K)=\cases{
0, &\quad $K \leq1/2$,\cr
(2K-1)/4, &\quad $1/2 \leq K \leq1$,\cr
K^2/4, &\quad $1 \leq K \leq2$,\cr
K-1, &\quad $K\geq2$.
}
\]
$P'$ then has a jump at $K=1/2$.

We have $\varphi(x) = \infty$ for $x < 1/2$ and $\varphi(x)=1/x$ for
$1/2\leq x \leq1$, and we assume that this formula also applies on
$[1,\infty)$. Then, $s(x)=2 \ln x$ for $x>1/2$, and
\[
g(y)=\cases{
e^{y/2}, &\quad $y>- 2 \ln2$,\cr
1/2, &\quad $y\leq-2\ln2$.
}
\]
Then, $\nu(dy)=dy/(8r)$ for
$y>- 2 \ln2$, $\nu(\{- 2 \ln2\})=1/(4r)$ and $\nu(dy)=0$ for $y<- 2
\ln2$. Consequently, the time change $\Gamma_u=\frac{1}{8r}O^+_u+
\frac{1}{4r}L_u^{-2 \ln2}$ is a linear combination of $O^+_u$ and
$L_u^{-2 \ln2}$, where $O^+_u$ is the amount of time spent by the
Brownian motion above $s(1/2)=-2 \ln2$ before time $u$.

It follows that for $x \geq1/2$, $\eta(x)^2 \equiv(x \sigma(x))^2 =
2rx^2$. As before we have
%
\begin{equation}\label{8.5b'example}
dX_t = rX_t\,dt + \sqrt{2r} X_t\,dB_t,\qquad t \leq H_{1/2}.
\end{equation}
It is easy to check using It\^{o}'s formula that $e^{-r \Gamma_u}
g(B_u)$ is a martingale in this case. The process $Z_t = B_{A_t}$ is
``sticky'' at $s(1/2)$ (this time in the sense of a~one-sided sticky
Brownian motion; see Warren~\cite{War}) and this property is inherited
by $X=s(Z)$.
\end{eg}

There is a third case, where $\varphi(\underline{x}+)< \infty$, but
$\varphi'(x+)=\infty$.

\begin{eg}
Suppose $\varphi(x) = 2 - \sqrt{2x-1}$ for $1/2 \leq x \leq5/2$ [and
$\varphi(x)=\infty$ for $x<1/2$]. Equivalently,
\[
P(K)=\cases{
0, &\quad $K \leq1/2$,\cr
2 - \sqrt{5-2K}, &\quad $1/2 \leq K \leq2$,\cr
K-1, &\quad $K\geq2$.
}
\]

For $1/2 < x < 5/2,$ we then have
\[
\eta(x)^2 = (x \sigma(x))^2 = 2r(2x-1) \bigl( 2 \sqrt{2x-1} + 1-x \bigr).
\]
It follows that although $X_t$ can hit $1/2$, the volatility at this
level is zero, and the drift alone is sufficient to keep $X_t \geq
1/2$.



\end{eg}

\subsection{\texorpdfstring{Kink in $P$ at $K^*$: $K^*<\infty$ and $P'(K^*-) < 1$}{Kink in P at K^*: K^* < infinity and P'(K^*-) < 1}}

In this case $\varphi'(x)$ is constant on an
interval $(\hat{x},x_0)$. This case is analogous to the one discussed
in Section~\ref{kink}.

\section{Extensions}

\subsection{No options exercised immediately}\label{noee}

In Hypothesis~\ref{hyp2}, in addition to~(i) and (ii), which are
enforceable by no-arbitrage considerations, we also assumed~(iii) that
there exists a finite strike $K^*$ such that for all strikes $K\geq
K^*$ the put option is exercised immediately. Since $K^*<\infty$ is
equivalent to $\varphi'(x_0)<0,$ it is apparent from the expression in
(\ref{phi-derivative}) that provided $\sigma$ is finite on some
interval $(x_1,x_2)$ where $x_0<x_1<x_2$ or, equivalently, $\nu$ gives
mass to some interval $(y_1,y_2)$ where $0<y_1<y_2$, this property
will hold. However, it is interesting to consider what happens when
this fails.

Suppose that $P(K)>K-x_0$ for all $K$ and that $\lim_{K\to\infty} P(K)
- (K-x_0)=0$. Then, $\varphi'(x_0)=0$, but $\varphi$ is strictly
decreasing on $(\underline{x},x_0)$. The measure $\nu$ places no mass
on $(0,\infty)$, the process $Z_t$ spends no time on $(0,\infty)$ and
$X_t$ never takes values above $x_0$. In particular, $X_t$ is reflected
(downward) at~$x_0$. The resulting model is consistent with observed
option prices, but not with the assumption that the discounted price
process is a (local) martingale. However, by allowing nonzero dividend
rates, we can find a model for which the ex-dividend price process is a
martingale and for which the model prices are given by $P(K)$; see
Section~\ref{zeror} below.

Now, suppose $\lim_{K} P(K) - (K-x_0)=\delta>0$. If $P(x_0)=x_0,$ then
$P(K)\equiv K$ and we have an extreme example which falls into this
setting. For $P$ as specified above we have that $\varphi(x)=1$ on
$(x_0-\delta,\infty)$. The measure $\nu$ places no mass on
$(-\delta,\infty)$ and $s(x)= x-x_0$ on this region. Except for time 0,
the process $Z_t$ spends no time in $(-\delta,\infty)$ and $X_t$ jumps
instantly to $x_0-\delta$, and thereafter spends no time above this
point. Alternatively, if $x_0$ is not specified, then this case can be
reduced to the previous case by assuming $x_0=K- \lim_K P(K)$.

\subsection{Nonzero dividend processes}
Until now, we have assumed that dividend rates are zero. However, if
dividend rates are a prespecified function of the asset price, then our
method adapts in a straightforward manner.

Given put prices $P(K),$ we recover $\varphi$ exactly as before from
the representation~\eqref{phi}. The unknown volatility $\sigma$ and
$\varphi$ are then related via the modified version of~\eqref{ODE}:
%
\begin{equation}\label{ode1}
\tfrac{1}{2}\sigma(x)^2 x^2 \varphi_{xx}+\bigl(rx -
q(x)\bigr)\varphi_x-r \varphi=0,
\end{equation}
where $q$ denotes the dividend rate. [We are assuming that under the
pricing measure, $X$ is governed by the stochastic differential
equation $dX_t = (r X_t - q(X_t))\,dt + \sigma(X_t)X_t\,dB_t$, where
$q$ is a known function.] The candidate $\sigma$ is then given by
%
\begin{equation}\label{sigmadef1}
\sigma^2(x) =
2\frac{ r \varphi - (rx - q(x)) \varphi_x}{x^2 \varphi_{xx} },
\end{equation}
at least where this quantity exists.\vadjust{\goodbreak}

Since $\varphi$ is convex by construction, a~necessary condition for
the prices $P(K)$ to be consistent with some model with dividend rate
$q$ is that $r \varphi- (rx - q(x))\varphi_x \geq0$. Then, in smooth
cases, where the existence of the diffusion with volatility $\sigma$
can be guaranteed, the analysis is complete. However, if $\varphi$ is
not strictly convex and twice differentiable, then some care may be
needed to define the diffusion associated with the candidate $\sigma$
given in~\eqref{sigmadef1}.

In keeping with our analysis in the previous sections, the most natural
approach for defining the (potentially generalized) one-dimensional
diffusion $X$ is via scale and speed. Note that if $\sigma$ is
sufficiently regular and $L_\sigma$ is the operator
\[
L_\sigma u =
\tfrac{1}{2} \sigma(x)^2 x^2 u_{xx} + \bigl(rx-q(x)\bigr) u_x - ru,
\]
then $L_\sigma\varphi= 0$. Moreover, we can find a second linearly
independent solution $\psi$ of $L_\sigma u = 0$ by the ansatz $\psi=
\varphi v$. This leads to the ODE
\[
\tfrac{1}{2} \sigma(x)^2 x^2 v_{xx} \varphi+ \sigma(x)^2 x^2 v_{x}
\varphi_{x} + \bigl(rx-q(x)\bigr) \varphi v_x = 0,
\]
which gives the unknown $v_x$ and its derivative in terms of $\varphi$
and $\sigma,$ and which has solution
%
\begin{eqnarray}\label{qscale}
v_x &= &\frac{A}{\varphi^2(x)} \exp\biggl( -
\int_{x_0}^x \frac{2(rz - q(z))}{z^2 \sigma(z)^2}\,dz \biggr)\nonumber\\ [-8pt]\\ [-8pt]
& = & \frac{A}{\varphi^2(x)} \exp\biggl( - \int_{x_0}^x \frac{
\varphi_{zz}(rz - q(z))}{r \varphi- (rz - q(z)) \varphi_z}\,dz \biggr).\nonumber
\end{eqnarray}
Note that the last expression is in terms of the dual function
$\varphi$ and does not involve $\sigma$ directly.

It is easily checked that the derivative of the scale function is given
by the Wronskian, so $s'(x) = \varphi\psi_x - \varphi_x \psi=
\varphi^2 v_x$. As before, the scale function can be used to determine
the inverse scale function $g$ and measure $\nu$. In turn, $\nu$ can be
used to determine the time change $\Gamma,$ and $X$ is given by the
formula $X_t = s(B_{A_t}),$ where $A$ is inverse to $\Gamma$. Thus, in
principle, the methods of this article extend directly to the case with
dividends, even in the irregular case [although further work is
necessary if $P$ is not strictly convex below $K^*$, whence $\varphi_x$
is not continuous, and the integral in~\eqref{qscale} is not well
defined]. However, we will not complete the analysis in this case and
instead will just make a remark and give a couple of examples.

\begin{remark*}
Whereas when dividends are zero
we have (e.g., from
Lem\-ma~\ref{varphiconvex}) that $\varphi$ is convex, this is not always
the case when dividends are positive. This means that the duality
between $P$ and $\varphi$ is more subtle. A~convex $P$ will lead to a
convex $\varphi$ and thence to a model which is consistent with the
perpetual put prices $P(K)$. However, starting with a model for which
$\varphi$ is not convex, we can still derive option prices $P$ from
expressions such as~\eqref{Ptilde}, but if we\vadjust{\goodbreak} now try to recover the
model from those prices, the duality lemma will lead to a function
$\tilde{\varphi} \neq\varphi$. Expressed differently, in the case with
dividends it is possible to have many time-homogeneous diffusion models
for which put prices are identical.
\end{remark*}

\begin{eg}
Suppose dividends are proportional so that $q(x)= \barq x$ with
$\barq\leq r$. Suppose further that $X_0=1,$ and $P(K)$ is given by
%
\begin{equation}
\hat P(K)=\cases{
\displaystyle\frac{K}{\beta+1}\bigl(\beta
K/(\beta+1)\bigr)^\beta, &\quad if  $K<K^*$,\vspace*{2pt} \cr
K-1, &\quad if  $K\geq K^*$,
}
\end{equation}
where $K^*=(\beta+1)/\beta$ for some positive $\beta$.

Then, $\varphi(x)=x^{-\beta}$. It follows that $\sigma^2(x) = 2(r +
(r-\barq) \beta)/(\beta(\beta+ 1))$. In this case it is clear that
$X$ is exponential Brownian motion and it is not necessary to calculate
the scale function. However, a~scale function can easily be computed
and is given by $s(x) = x^{c} - 1$, where $c = (r -
(r-\barq)\beta^2)/(r + (r-\barq) \beta)$.
\end{eg}

\begin{eg}
Suppose that, as in Example~\ref{eg:kink}, $X_0=1$ and $P(K)$ is
given by
\[
P(K) = \cases{
{K^2}/{8,}
&\quad $0 < K \leq27/32$,\vspace*{1pt} \cr
4K^3/27, &\quad $27/32 \leq K \leq3/2$,\vspace*{1pt} \cr
(K-1), &\quad $3/2 \leq K$.
}
\]
This time, however, we assume that there are proportional dividends
with constant of proportionality $\barq$ (with $\barq<r$). As in
Example~\ref{eg:kink}, we find that (with $\lambda= 4/3$)
\[
\varphi(x) = \cases{
{2}x^{-1}, &\quad $0 < x \leq 27/64 = \lambda^{-3}$,\vspace*{1pt} \cr
4 \lambda^3 - 2 \lambda^6 x, &\quad $27/64= \lambda^{-3} < x \leq 9/16 = \lambda^{-2}$,\vspace*{1pt} \cr
x^{-2}, &\quad $x> 9/16 = \lambda^{-2}$.
}
\]
Then, from $s' = \varphi^2 v'$ and~\eqref{qscale} we find that $s$ is
linear over $I=[27/64,9/16]$ and, more generally, a~choice of $s$ can
be obtained by integrating
\[
s'(x) = \cases{
\lambda^b x^{-2(r-\barq)/(2r-\barq)},
&\quad $0 < x \leq27/64$,\vspace*{1pt} \cr
\lambda^{2c}, &\quad $27/64 < x \leq9/16$,\vspace*{1pt} \cr
x^{-c}, &\quad $x>9/16$,
}
\]
where now $c= 6(r - \barq)/(3r-2 \barq)$ and $b= 6r(r - \barq)/[(3r-2
\barq)(2r - \barq)]$.
\end{eg}

\subsection{Time-homogeneous processes with known volatility and
unknown interest rate or unknown dividend processes} \label{zeror}

In the main body of the paper we have assumed that the interest rate
$r$ is a given positive constant, that dividend rates are zero and that
$\sigma$ is a function to be determined. In the last section we
generalized this analysis to allow for a known, nonzero, dividend rate.
We will now argue that the same ideas can be used to find other
time-homogeneous models consistent with observed perpetual put\vadjust{\goodbreak} prices,
whereby the volatility function is given, and either a state-dependent
dividend rate or a state-dependent interest rate is inferred.

Suppose $X$ has dynamics $dX_t = (r(X_t)X_t - q(X_t))\,dt +
\sigma(X_t)X_t\,dB_t$. Given put prices $P(K)$ as before, define
$\varphi$ via $\varphi(z) = \inf_{K:K\geq z} (K-z)/P(K)$. The
relationship between $\varphi$ and the characteristics of the price
process $X$ are then such that $\varphi$ solves $L\varphi= 0$, where
$L$ is given by
\[
Lu = \tfrac{1}{2} \sigma(x)^2 x^2 u_{xx} + \bigl(xr(x)-q(x)\bigr) u_x - r(x)u.
\]
Note that we now allow for any of $\sigma$, $q$ or $r$ to depend on
$x$. Until now we have assumed that $r$ is constant and $q$ is zero
(except in the last section, where $q$ was known but nonzero) and
solved for $\sigma$, but we can alternatively assume that $\sigma(x)$
is a given function and $r$ is a positive constant, and solve for $q$,
or assume that $q$ and~$\sigma$ are given, and solve for $r$.

For example, if $r$ and $\sigma$ are given constants, then the dividend
rate process is given by
\[
q(x) = x r + \frac{ x^2 \sigma^2 \varphi_{xx} - 2 r \varphi
}{ 2 \varphi_x} .
\]
If $q$ is negative, this should be thought of as a
convenience yield.

By allowing for dividend processes which are singular with respect to
calendar time and which are instead related to the local time of $X$ at
level $x_0$, it is possible to construct candidate price processes
which spend no time above $x_0$. For example, if $\tilde{L}$ is the
local time at $1$ of $X$, and if
\[
\frac{dX_t}{X_t} = dB_t + r\,dt - \frac{d\tilde{L}_t}{2},
\]
then $X_t$ reflects at 1, and if $\hat{\varphi}(x) =
(\E^1[e^{-rH_x}])^{-1}$ for $x<1$, then $\hat{\varphi}^\prime(1-)=0$.
This gives an example of a model consistent with the class of option
prices described in Section~\ref{noee}.

\subsection{Recovering the model from perpetual calls}

The perpetual American call price function $C\dvtx[0,\infty)\rightarrow
[0,x_0]$ must be nonincreasing and convex as a function of the strike
$K$, and must satisfy the no-arbitrage bounds $(x_0-K)^+ \leq C(K) \leq
x_0$.

If there are no dividends (and if $e^{-rt}X_t$ is a martingale), then
the perpetual call prices are given by the trivial function $C(K)=x_0$.

So, suppose instead that the (proportional) dividend rate $\barq$ is
positive. Let $\hat{\psi}$ be the increasing positive solution to
\[
\tfrac{1}{2} x^2 \sigma(x)^2 \hat\psi'' + (r-\barq)x \hat\psi' - r
\hat\psi= 0,
\]
normalized so that $\hat{\psi}(x_0)=1$. Then, for
$z>x$, $\E^x[e^{-r H_z}] = \hat{\psi}(x)/\hat{\psi}(z)$, and call
prices in a model where $dX_t = (r-\barq)X_t\,dt +
\sigma(X_t)X_t\,dB_t$ are given by
\[
\hat{C}(K) = \sup_{\tau} \E^{x_0}[e^{-r \tau}(X_\tau- K)^+] =
\sup_{x :x\geq x_0 } \frac{(x-K)}{\hat{\psi}(x)} .
\]

\begin{eg*}
Suppose $X$ solves $(dX_t/X_t) =
(r-\barq)\,dt + \sigma\,dB_t$ with
$X_0=x_0$. Then, $\hat\psi(x) = (x/x_0)^{\gamma},$ where $\gamma=
\beta_+$ and
\[
\beta_{\pm} = - \biggl( \frac{r-\barq}{\sigma^2} - \frac{1}{2}
\biggr)
\pm\sqrt{ \biggl( \frac{r-\barq}{\sigma^2}
- \frac{1}{2} \biggr)^2+\frac{2r}{\sigma^2} }.
\]
Note that since $\barq>0$ we have $\gamma>1$. Note also that
$\hat{\varphi}(x) = (x/x_0)^{\beta_-}$.

The corresponding call prices are given by
\[
\hat C(K) = x_0^{\gamma} \sup_{x:x \geq x_0 } \{ {(x-K)}x^{-\gamma}
\},
\]
which for $K \leq(\gamma-1)x_0/\gamma$ gives $\hat
C(K)=(x_0-K)$, and for $K > (\gamma-1)x_0/\gamma$ gives
\[
\hat C(K) = x_0^\gamma
{\gamma^{-\gamma}}(\gamma-1)^{\gamma-1}K^{1-\gamma}.
\]
\end{eg*}

The example discusses the forward problem, but the discussion of the
inverse problem is similar to that in the put case. Given perpetual
call prices $C(K)$, for $x>x_0$ we can define $\psi$ via $\psi(x)=
\inf_{K: K\leq x} (x-K)/C(K)$ and then construct a~triple $\sigma(x),
q(x), r(x)$ so that
\[
\tfrac{1}{2} x^2 \sigma(x)^2 \psi'' + \bigl(xr(x)-q(x)\bigr) \psi' - r(x) \psi
= 0 .
\]
By combining information from put and call prices, it is possible to
determine a candidate model which simultaneously matches both puts and
calls. The information contained in the perpetual puts determines the
volatility below $x_0,$ and the information contained in the perpetual
calls determines the volatility above $x_0$. However, for this
candidate model to return the put and call prices, there is an
additional consistency condition. For a discussion of this condition in
the smooth case, see Alfonsi and Jourdain~\cite{AJ08b}, Proposition
4.6.

\begin{appendix}

\section*{Appendix: Proofs}\label{app}

\subsection{Duality}\label{app1}
\mbox{}
\begin{pf*}{Proof of Lemma~\ref{dualitylemma}} It is clear that $g$
is nonnegative and nondecreasing since $f$ is positive and
nonincreasing. The lower bound on $g$ follows from choosing
$z=z_0\wedge k$ in~\eqref{g}, and the upper bound follows since $f$ is
nonincreasing. To show that $g$ is convex, first note that $g(k)$ is
minus the reciprocal of the slope of the tangent of the function $f$
which passes through the point~$(k,0)$.

For two given points $k_1$ and $k_2$ with $k_1<k_2$, let $l_1(z)$ and
$l_2(z)$ be the corresponding tangent lines. Let $k= \lambda
k_1+(1-\lambda)k_2$ for some $\lambda\in(0,1)$, and let~$l(z)$ be the
line through the point $(0,k)$ and the intersection point of $l_1$ and
$l_2$ (cf. Figure~\ref{fig3}).
%
\begin{figure}

\includegraphics{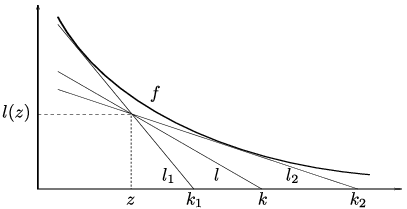}

\caption{The lines $l_1$, $l_2$ and $l$.}\label{fig3}
\end{figure}
If the intersection\vadjust{\goodbreak} point is denoted $(z,l(z))$, then the convexity of
$f$ guarantees that
\[
g(k)\leq \frac{k-z}{l(z)}=\frac{(1-\lambda) k_1-(1-\lambda) z}{l_1(z)}
+\frac{\lambda k_2-\lambda z}{l_2(z)}=(1-\lambda) g(k_1)+
\lambda g(k_2),
\]
which proves that $g$ is convex.

To prove the self-duality, let $z\leq z_0$. By the definition of $g,$
we have that $g(k)\geq(k-z)/f(z)$ for all $k\geq z$. Consequently,
\[
F(z)=\sup_{k\geq z}\frac{k-z}{g(k)}\leq f(z).
\]
For the reverse inequality, let $z\leq z_0$ and let $l$ be a tangent
line to $f$ through the point $(z,f(z))$ (such a tangent is not
necessarily unique if $f$ has a kink at~$z$). Assume that the point
where $l$ intersects the $z$-axis is given by $(k^\prime,0)$. Then,
$g(k^\prime)=(k^\prime-z)/f(z)$, so
\[
F(z)=\sup_{k\geq z}\frac{k-z}{g(k)}\geq
\frac{k^\prime-z}{g(k^\prime)}=f(z),
\]
which completes the proof of
(ii). The proof of (b) can be constructed along the same lines.
\end{pf*}

\subsection{Time changes of local martingales}\label{app2}

\begin{proposition}\label{tc-m}
 Suppose $(\gamma_u)_{u \geq0}$ is a martingale with
respect to the filtration $\mathbb G=({\mathcal G}_u)_{u \geq0}$, and
$A_t$ is an increasing process such that $A_t$ is a finite stopping
time with respect to $\mathbb G$ for each $t$. Define $\tilde{\gamma}_t
= \gamma_{A_t}$ and $\tilde{\mathcal G}_t = {\mathcal G}_{A_t}$. In
general $(\tilde{\gamma_t})_{t \geq0}$ is not a martingale. However,
if $\gamma$ is a bounded martingale, then $\tilde{\gamma}$ is a bounded
martingale.
\end{proposition}

\begin{pf}
Given a Brownian motion $B$, for $b>0$, let $H^B_b$ be the first
hitting time of level $b$. Then, $(\tilde{B}_b)_{b \geq0}$ defined via
$\tilde{B}_b \equiv B_{H^B_b}$ is not a martingale.

However, if $\gamma$ is bounded, then
$\E[\tilde{\gamma}_t|\tilde{\mathcal G}_s] = \E[\gamma
_{A_t}|{\mathcal
G}_{A_s}] = \gamma_{A_s} = \tilde{\gamma}_s$, by optional sampling.\vadjust{\goodbreak}
\end{pf}

Suppose now that we are in the setting of Section~\ref{sec-timechange},
where $Z_t$ is constructed from the Brownian motion $B$. In particular,
$\Gamma_u$ is an increasing additive functional of $B$, and $A$ is the
right-continuous inverse to $\Gamma$.

\begin{pf*}{Proof of Lemma~\ref{tc-lm}} Intuitively, a~time change
of Brownian motion is a local martingale, but if the additive
functional $\Gamma$ is constant when $B$ is in $[a,\infty),$ then the
resulting process spends no time above $a$ and reflects there. To
maintain the local martingale property we need either that the
time-changed process never gets to $a$, or
that there are 
arbitrarily large values at which $\Gamma$ is strictly increasing.

If $[\underline{z}_\nu,\overline{z}_\nu]$ is a bounded interval, then
$A_\infty\leq H^B_{\underline{z}_\nu} \wedge H^B_{\overline{z}_\nu}$
and $(Z_t)_{0 \leq t < \infty} = (B_{A_t})_{0 \leq A_t < A_\infty}$ is
a bounded martingale, by Proposition~\ref{tc-m}.

Now, suppose $(\underline{z}_\nu,\overline{z}_\nu) = \R$ and suppose
that for each $a$, $\nu$ assigns mass to every set $(a,\infty)$ and
$(-\infty,-a)$.

We have $A_{\Gamma_t} \geq t$ with equality when $\Gamma$ is strictly
increasing at $t$. Let $\{ a^+_n\}$ and $\{ a^-_n\}$ be two sequences
converging to $+\infty$ and $-\infty$, respectively, so that $\nu$
assigns mass to any neighborhood of $a^+_n$ and $a^-_n$, and set $H_n =
\inf\{u\dvtx B_u \notin(a^-_n,a^+_n) \}$. Then, $\Gamma$ is strictly
increasing at $H_n$. Set $T_n = \Gamma_{H_n}$. Then, $A_{T_n} = H_n$.
Note that $\Gamma_u$ increases to infinity almost surely, and hence
$\Gamma_{H_n} \uparrow\infty$. Under our hypothesis, $(Z^{T_n}_t)_{t
\geq0}$ given by
\[
Z^{T_n}_t := Z_{t\wedge T_n}= B_{A_{t \wedge T_n}} = B_{A_t \wedge H_n}
\]
is a bounded martingale. Hence, $T_n$ is a localization sequence for
$Z$.

The mixed case can be treated similarly.
\end{pf*}

\begin{pf*}{Proof of Lemma~\ref{lem-lm}} For $y\in(s(0),s(\infty))$,
set $H(y)= \ln(g(y)/x_0)$ and write $h(y) = H'(y) = g'(y)/g(y)$. If
$g$ is not differentiable at $y,$ then we take the right derivative,
which exists since $g$ is convex. (We use a similar convention for~$f$,
$h$ and $j$ defined below.) Then,
\[
\nu(dy) = \frac{1}{2r} \frac{g''(dy)}{g(y)} =
\frac{1}{2r} \bigl(h'(dy) + h(y)^2 \,dy \bigr)
\]
and, as usual,
$\nu(\{y\}) = \infty$ for $y\notin[s(0),s(\infty)]$. Note that in the
case where $g$ is not twice differentiable, we have $H''(dy) \equiv
h'(dy) = g''(dy)/g(y) - h(y)^2\,dy$ so that $H''$ exists in a
distributional sense.

We have $H(y) = \int_0^y h(v)\,dv = \ln(g(y)/x_0)$ and
\[
\Gamma_u = \frac{1}{2r} \int_{\R} L^y_u \bigl( H''(dy) +
H'(y)^2\, dy \bigr) .
\]
Let $\xi$ be the first explosion time of $\Gamma$.
Then, by the It\^{o}--Tanaka formula (e.g., Revuz and Yor~\cite{RY}, Theorem
VI.1.5), for $u<\xi$,
\begin{eqnarray*}
H(B_u) & = & \int_{0}^{u} H'(B_s)\,dB_s + \frac{1}{2} \int_{\mathbb R}
L^y_{u}H''(dy) \\
& = & \int_{0}^{u}h(B_s)\,dB_s - \frac{1}{2} \int_{\mathbb R} L^y_{u}
h(y)^2\,dy + r \Gamma_{u}.
\end{eqnarray*}
Thus, $g(B_{u} ) = x_0 e^{H(B_{u})} = {\mathcal E}(h(B) \cdot B)_{u}
e^{r \Gamma_{u}}$, where $\mathcal E$ denotes the Dol\'{e}ans
exponential, and $e^{-r \Gamma_{u}} g(B_{u})$ is a local martingale. It
follows that $M$ is a local martingale.

Now, define $J(y) = \int_0^y j(v)\,dv = \ln f(y)$ and
\[
\tilde{\Gamma}_u = \frac{1}{2r} \int_{\R} L^y_u \bigl( J''(dy) +
J'(y)^2\,dy \bigr) .
\]
Again, $J''(dy)= f''(dy)/f(y) - j(y)^2\,dy$
exists in the distributional sense, even if $j(y)$ is not continuous.
By exactly the same argument as above, we find that $f(B_{u}) =
e^{J(B_{u})} = {\mathcal E}(j(B) \cdot B)_{u} e^{r \tilde{\Gamma}_{u}}$
and $e^{-r\tilde{\Gamma}_{u}} f(B_u)$ is a local martingale.

It remains to show that $\Gamma_u = \tilde{\Gamma}_u$. Define $L(y) =
(f(y)g(y))^{-1}$ so that $L$ is continuous and right-differentiable.
[We write $L'(y)$ for this right-derivative when the derivative is not
well defined.] Then, $L'(y)/L(y) = - g'(y)/g(y) - f'(y)/f(y) = - (H'(y)
+ J'(y))$ and
\begin{eqnarray*}
J'(y) - H'(y) &=& \frac{\varphi'(g(y))g'(y)}{\varphi(g(y))} -
\frac{g'(y)}{g(y)} = \frac{g'(y) [ g(y) \varphi'(g(y)) -
\varphi(g(y))]}{g(y)\varphi(g(y))}\\
&=& - \frac{g'(y)s'(g(y))}{g(y)f(y)} = - L(y).
\end{eqnarray*}
We have that $J'(y) - H'(y)$ is (right-) differentiable, even if
separately $J'$ and $H'$ are not, and
\[
\bigl(J'(y) - H'(y)\bigr)' = \bigl({H'(y)+J'(y)}\bigr)({L(y)}) =
H'(y)^2 - J'(y)^2 .
\]
Finally, since $L^y_u$ is a bounded continuous
function with compact support for each fixed $u$, we conclude that
$\Gamma_u = \tilde{\Gamma}_u$.
\end{pf*}

\begin{pf*}{Proof of Corollary~\ref{e-rtX=lm}} Recall that in our
setting, $\Gamma$ defined via~\eqref{Gammadef} grows without bound and
is continuous, at least until $B$ hits $s(0)$ or $s(\infty)$. Thus, if
$\xi$ denotes the first explosion time of $\Gamma$, then the inverse
function $A$ is defined for every $t,$ and $A_t = \xi$ for $t \geq
{\Gamma_\xi}$. Then, using the extension of the definition of
$\tilde{M}$ beyond $\Gamma_\xi$ as necessary, we have
\[
\tilde{M}_t = e^{-rt} X_t =
\cases{
M_{A_t}, &\quad $t \leq\Gamma_\xi$, \cr
M_{\xi}, &\quad $t > \Gamma_\xi$.
}
\]
%

Recall that $\varphi$ is extended to $(x_0,\infty)$ in such a way that
$\lim_{x \uparrow\infty} \varphi(x)=0$. Therefore, either
$s(\infty)<\infty$ and $\nu$ assigns infinite mass to all points $z >
s(\infty) = \overline{z}_{\nu}$, or $s(\infty)=\infty$ and there exists
a sequence $a_n \uparrow\infty$ such that $\nu$ assigns mass to any
neighborhood of $a_n$.

Suppose that the second case obtains. If $s(0)>-\infty,$ then
$H_{s(0)}^B=\xi<\infty$, otherwise $\xi=\infty$. On
$H_{a_n}^B<H^B_{s(0)}=\xi$, $\Gamma_u$ is strictly increasing at $u =
H_{a_n}^B$ and $A_{\Gamma_{H_{a_n}^B}}= H_{a_n}^B$. Set
%
\begin{equation}\label{Tndef}
 T_n = \cases{
  \Gamma_{H_{a_n}^B},
&\quad $H_{a_n}^B < H^B_{s(0)}$,\vspace*{2pt} \cr
\infty, &\quad $H_{a_n}^B > H^B_{s(0)}$,
}
\end{equation}
where the second line is redundant\vspace*{-2pt} if $s(0)=-\infty$. Then, $A_{T_n} =
H^B_{a_n} \wedge\xi$ is such that $\tilde{M}^{T_n}_t :=
\tilde{M}_{t\wedge T_n} =
M_{A_t
\wedge\xi\wedge H^B_{a_n}} \leq g(a_n)$ and $T_n$ is
a reducing sequence for~$\tilde{M}$.

Now, suppose $s(\infty)<\infty$ and $\bar{g} = \infty$. Choose $a_n
\uparrow s(\infty)$ such that $\nu$ assigns mass to any neighborhood of
$a_n$. Then, on $H^B_{s(\infty)}<H^B_{s(0)}$, we have,\vspace*{1pt} by the argument
after Lemma~\ref{lem-lm}, that $\Gamma_{H^B_{a_n}} \uparrow\infty$
almost surely, and the argument proceeds as before with $T_n$ given by
(\ref{Tndef}) being a reducing sequence.

Finally, suppose $s(\infty)<\infty$ and $\bar{g} <\infty$. Then, $M$
is bounded by $\bar{g}$ and $\tilde{M}$ is a martingale.
\end{pf*}
\end{appendix}

\section*{Acknowledgments}
We thank Eberhard Mayerhofer and an anonymous referee
for their careful reading.

%

\printaddresses

\end{document}